%% file: main_v6_arxiv_final.tex
\documentclass[english,11pt,a4paper]{article}

\usepackage[
    backend=biber,
    style=alphabetic,
    maxcitenames=100,
    maxbibnames=100,
    backref=false
]{biblatex}
\usepackage[margin=22mm]{geometry}

\usepackage{amsmath,amssymb,amsthm}
\usepackage{thmtools}
\usepackage{subcaption}
\usepackage{array,blkarray,longtable,multirow}
\usepackage{booktabs} 
\usepackage{adjustbox}
\usepackage{enumitem}

\makeatletter
\let\original@algocf@latexcaption\algocf@latexcaption
\long\def\algocf@latexcaption#1[#2]{%
	\@ifundefined{NR@gettitle}{%
		\def\@currentlabelname{#2}%
	}{%
		\NR@gettitle{#2}%
	}%
	\original@algocf@latexcaption{#1}[{#2}]%
}
\makeatother

\newcommand{\piv}[1]{\mathrm{piv}\left(#1\right)}
\newcommand{\before}{\triangleleft}
\newcommand{\after}{\triangleright}
\newcommand{\col}[2]{#1\left[#2\right]}

\newcommand{\str}[1]{\mathrm{st}\left(#1\right)}

\newcommand{\nameshort}{\mbox{SiRUP}}

\newcommand{\deyhoualgo}{\mbox{\textsc{FastZigzag}}}
\newcommand{\phatv}{\mbox{\textsc{PHAT-S}i\textsc{RUP}}}

\usepackage[algo2e,vlined,ruled,linesnumbered]{algorithm2e} 
\SetKwIF{If}{ElseIf}{Else}{if}{}{else if}{else}{end if}%
\SetKwFor{While}{while}{}{end while}%
\SetKwFor{For}{for}{}{end for}%
\DontPrintSemicolon

\newcolumntype{C}{>{\centering\arraybackslash$} p{.3cm} <{$}}

\makeatletter
\renewcommand{\nllabel}[1]
{{\let\@currentlabel\algocf@currentlabel
		\let\@currentcounter\algocf@currentcounter
		\label{#1}}}%

\renewcommand{\algocf@nl@sethref}[1]{%
	\renewcommand{\theHAlgoLine}{\thealgocfproc.#1}%
	\hyper@refstepcounter{AlgoLine}%
	\gdef\algocf@currentlabel{#1}%
	\gdef\algocf@currentcounter{AlgoLine}%
}%
\makeatother

\usepackage[hidelinks]{hyperref}
\usepackage[noabbrev,nameinlink,capitalise]{cleveref}

\theoremstyle{plain}
\newtheorem{theorem}{Theorem}[section]
\newtheorem{lemma}[theorem]{Lemma}
\newtheorem{proposition}[theorem]{Proposition}

\theoremstyle{definition}
\newtheorem{remark}[theorem]{Remark}
\newtheorem{example}[theorem]{Example}

\makeatletter
\let\original@algocf@latexcaption\algocf@latexcaption
\long\def\algocf@latexcaption#1[#2]{%
	\@ifundefined{NR@gettitle}{%
		\def\@currentlabelname{#2}%
	}{%
		\NR@gettitle{#2}%
	}%
	\original@algocf@latexcaption{#1}[{#2}]%
}
\makeatother

\DeclareMathOperator{\rank}{rank}

\newcommand\define[1]{\textbf{#1}}

\usepackage[dvipsnames,table]{xcolor}
\usepackage{array,colonequals,tikz,tikz-cd,fp}
\usetikzlibrary{calc, arrows.meta, arrows, shapes.geometric}
\usepackage{mathtools} 
\usepackage{extarrows}

\definecolor{color4}{HTML}{9467bd}    
\definecolor{color5}{HTML}{8c564b}    

\definecolor{dimzerocolor}{HTML}{3FC7F4}
\definecolor{dimonecolor}{HTML}{F55881}
\definecolor{dimtwocolor}{HTML}{F5E727}
\definecolor{dimthreecolor}{HTML}{1F43C2}

\definecolor{color0}{HTML}{4285f4}    
\definecolor{color1}{HTML}{0f9d58}    
\definecolor{color2}{HTML}{f4b400}    
\definecolor{color3}{HTML}{db4437}    
\definecolor{color7}{HTML}{737373}    

\usepackage[textsize=scriptsize]{todonotes}

\begin{document}

\title{Pruning vineyards: updating barcodes and representative cycles by removing simplices} 
\author{Barbara Giunti\footnote{\href{mailto:bgiunti@albany.edu}{\texttt{bgiunti@albany.edu}} SUNY - University at Albany, NY, USA, \& Graz University of Technology, Austria} 
	\and J\=anis Lazovskis\footnote{\href{mailto:janis.lazovskis@lu.lv}{\texttt{janis.lazovskis@lu.lv}} University of Latvia, Riga, Latvia, \& University of Aberdeen, Scotland, United Kingdom}}
\date{\today}

\maketitle
	
\begin{abstract}
The barcode of a filtration and its representative cycles encode rich information often useful in data analysis.
However, obtaining them can be computationally expensive.
Therefore, it is useful to have methods that update them if the associated filtration undergoes small changes.
There are already efficient algorithms updating a barcode if simplices exchange entrance order or are added, but not if simplices are removed.
We provide an implementation to update a reduced boundary matrix when simplices in the filtration are removed.
Our algorithm, the Simplicial Removal Update Procedure (SiRUP), intrinsically updates also the representative cycles, and is compatible with the clearing optimizations.
We show that the complexity of our algorithm is lower than recomputing the barcode from scratch and that the number of executed matrix column additions is minimal, with both theoretical and experimental methods.
\end{abstract}



\section{Introduction}\label{sec_intro}

Topological data analysis (TDA) is a discipline at the intersection of mathematics, computer science, and data analysis, that uses tools from (traditionally pure) mathematics to analyze data.
Its main difference from other approaches in applied mathematics is that it is a features extractor: it extract shape-related information from the data, such as the presence (or absence) of connected components, loops, and cavities.
The mathematical name for this kind of information is \define{homology}.
Thus, more than correlating different characteristics of a dataset (as statistics would do), it directly describes properties of it.
The information extracted can be then interpreted directly or further studied with statistical or machine learning methods.
TDA has many advantages: it does not require the input to be smooth, it has virtually no limitation on which type of input can be processed, and the extracted information is typically \textit{stable}, i.e., from similar data will output similar information.
Since any data carries some noise, this latter advantage is particularly useful.
One subdiscipline of TDA is specially concerned with stability: persistence theory.
The idea behind persistence theory is to extract information (typically, homology) for a varying set of thresholds in the data, such as distances in a point cloud, pixel values in a binarized image, function values in a scalar field, and so on.
Thus, the information extracted is the evolution of a certain property in the data, which is usually more stable because it does not depend on a chosen value for the thresholds.
If the property of interest was homology, we call the computation of its evolution \define{persistent homology (PH)}.

Persistent homology is a tremendously successful tool in topological data analysis, with hundreds of practical applications \cite{donut}.
Arguably, one of the reasons for its success is the computability of its main invariant, the \define{barcode}, \cite{bauer2021ripser,phat_paper,roadmap}, which represents the evolving homology of a dataset as a collection of life spans of homological features, or \define{bars}.
Barcodes can be equivalently represented as \define{persistence diagrams}: 2-D plots whose x-axis represents the values at which a (homological) feature starts, and y-axis represents the values at which that feature ends.
Since they are mathematically equivalent, we will use the terms ``barcode'' and ``persistent diagram'' interchangeably in the text.
The \define{representative cycles} are collections of \define{simplices} (which are subsets of data points) that generate the topological features encoded in the barcode, and have been demonstrably useful in practice \cite{tight_cycles,cycles_userguide}.
For practical purposes, the worst-case complexity of the barcode computation is cubic in the number of simplices, but even considering the best theoretical worst-case complexity, which is matrix-multiplication time \cite{matrix_multi}, the the runtime complexity becomes problematic for very large datasets.
It is therefore useful to have methods to update an already computed barcode if the dataset changes a little, instead of recomputing it from scratch.

The barcode computation requires simplices to be ordered, or \define{filtered}, and in any ordering there are fundamentally three changes that may happen: simplices are added, simplices exchange order, or simplices are removed.
In the literature, insertion and order swap of simplices has already been studied \cite{vines_and_vineyards,luo2023accelerating}, but no explicit, systematic work has been done on updating the barcode and the representative cycles by removal of simplices.
The removal of simplices from a complex may be encoded by a zigzag filtration \cite{dey2022updating,dey2023computing}, but that is not the most efficient method to remove simplices (see \cref{ex_dey-hou-comparison}).
The approaches in \cite{vines_and_vineyards,luo2023accelerating} allow for removal only by swapping the simplices to the end of the filtration and then dropping them, which incurs more operations than are strictly necessary (see \cref{sec_comparison} for an extended discussion).
In this paper, we fill this gap by providing an efficient method to update a barcode and the representative cycles that generate it by removing simplices.

We provide an algorithm that takes as input the boundary matrix $D$ factored as $D=RV$, where $V$ is the \define{operations matrix}, and a list $L$ of simplices to be removed.
This factorization is the usual output of barcode algorithms, as the barcode and almost all the representative cycles\footnote{i.e., the ones corresponding to finite bars in the barcode} can be read directly from $R$.
Our algorithm returns $R'$ and $V'$, where $D'=R'V'$ is the matrix of the original filtration without the simplices in $L$.
We show that our method of computing $R'$ and $V'$ from $R$ and $V$ has better theoretical complexity than factoring $D'$ into its two factors.
Moreover, we provide an implementation and experiments showing that also in practice we can expect our method to be faster.

\subsection{Motivation and applications}
\label{sec_motivations_applications}

\subparagraph{Neuroscience.}
Computational neuroscience has recently been harnessed with topological tools \cite{nnneighbourhoods,morphomics,ccliques}, which allow for a dramatic simplification.
A common approach is to consider neurons as data points and synapses as (directed) connections between the data points; or a cell body and branchings as data points and the cell material in between as connections.
This describes a directed graph, whose \textit{clique complex} is a collection of simplices with a natural filtration by dimension and ordering.
Neuron activity (either firing or not firing) gives a sequence of filtrations over time, defined by the active subcomplex.
When a collection of neurons fires in one time interval and does not fire in the next, the change is reflected as simplices removed from the filtration.
Datasets in neuroscience often have several thousands of neurons, with current experiments often carried out on smaller subsets because the computational complexity is too high.
Therefore, it is useful to have the possibility to ``pay the price'' by computing once the whole the barcode of datasets, and then perform only small and local update when different neurons fire.

\subparagraph{Manifold learning.} 
As a machine learning method, manifold learning reduces the dimensionality of the data while preserving important features.
Non-linear dimensionality reduction is often performed by \textit{auto-encoders}, which can be paired with PH to preserve more information, in particular topological information \cite{bastianautoencoders}.
The simplest auto-encoders attempt to reduce the dimension of the manifold globally, without respect to local features.
Using multi-charts auto-encoders is a way to keep local information, taking care to appropriately relate the barcodes of each chart.
However, there is some fine-tuning to be done to choose the right number and type of charts, balancing efficiency and information.
Every time charts are modified, the barcode needs to be recomputed.
With our method, it is enough to make an initial computation which can be continuously updated whenever charts are modified.

\subparagraph{Streaming, dynamic, and noisy data.} 
In settings where data is incomplete, or regularly being updated, any sort of topological signature on the data will have to be updated as frequently as the underlying data is updated \cite{sliding_windows}.
This area has only recently been developing \cite{time_varying_reeb_graphs,facundo_dynamic},
but its general setting of continuously responding to new inputs makes it a fertile area for many directions, emphasized by the proliferation of data and data channels in the 21st century.
This can also be applied to situations where the full topological analysis of a dataset is not possible because the total size is too big.
In this case, one can begin to analyze a subset of the input, and then add and remove data to look for meaningful topological information without exceeding the computational capability.

\subsection{Related work}
\label{sec_relatedwork}
Changes in PH over a time parameter have been widely studied \cite{cerrifrosini,time_varying_reeb_graphs,hickok_bundle,time_varying_merge_trees}, with algorithmic advances for dynamic computations laid in \cite{vines_and_vineyards,facundo_dynamic}, introducing methods to update the reduced boundary matrix when the order of simplices changes.
Building on this work, the authors of \cite{luo2023accelerating} provide efficient methods to update the barcode of any filtration, if simplices are inserted, their ordered changed, or if the simplices at the end of the filtration are removed.
This is the only case we are aware of a work explicitly discussing the updating of barcode following the removal of simplices, though it uses the order changing method of \cite{vines_and_vineyards} to push the simplices to the end of the filtration and then drops them.
Indeed, the work of Cohen-Steiner--Edelsbrunner--Morozov \cite{vines_and_vineyards} underlies most implementations for computing persistent homology with specific methods for dynamic filtrations (for example, \texttt{pl-vineyard} of \cite{dionysus} and \texttt{remove\_maximal\_simplex} of \cite{maria2014gudhi}).
The differences in procedure and efficiency in comparison with our method are detailed in \cref{sec_vinescomparison}, and the implementation specifics are presented in \cref{sec_implementations}.

In \cite{luo2023accelerating}, removal is reduced to swapping: keeping $R$ upper triangular while moving the simplex to be removed to the end of the filtration implies that simply dropping the last column maintains the correct pivot pairing of the modified filtration, as the last column can not be added to any other columns.
While, as shown in \cite{luo2023accelerating}, it is possible to implement a removal working around the standard barcode algorithm (\cref{alg_sba}), this is not very efficient, as many operations are canceling each other.
Our method skips several of these redundant operations by considering the removal in its own right.

Recent work by Dey--Hou \cite{dey2022updating, dey2023computing} computes the barcode of a \textit{zigzag} filtration, generalizing from a linear filtration order.
Their method may be adapted to update the barcode by removal of simplices, by inserting at the end of the filtration backward maps, each map removing a simplex.
Their \textit{forward switch} and \textit{inward contraction} moves (see \cite{dey2022updating, dey2023computing}) allow their methods to be used across a wide range of filtrations.
However, their approach doubles the length of the filtration as it is concerned with a wide class of changes.
While both their and our methods require in the worst case an amount of column additions that is linear in the number of simplices, we demonstrate ours always performs less or at most the same number of column additions, and has a more direct computation (see \cref{sec_zigzagcomparison}).

\subsection{Complexity, pseudocodes, and implementation}
\label{subs_(pseudo)code}
The study of the efficiency of algorithms is layered, involving the theoretical behavior of the pseudocode and the practical behavior of the actual implementation, as well as the worst-case complexity and the average complexity.
Conclusions drawn from these different layers may not agree: the easiest conclusions to obtain (often theoretical behavior and worst-case complexity) may be far from what happens in real-case scenarios.
For example, most of the algorithms computing the barcode have cubic worst-case complexity in the number of simplices (where the number of simplices is often a power of the number of input points), but in practice they are almost linear \cite{ripser_software,phat_paper,roadmap}.

In this work, we strive to obtain a bigger picture by studying many of these aspects: first, we analyze the theoretical behavior of the pseudocodes and their worst-case complexities in \cref{sec_removing}, followed by a discussion and benchmarking of our implementation to better understand the practical behavior, in \cref{sec_experiments}.
As a result, when presenting pseudocodes, we do not consider some coding aspects that can heavily affect the complexity, such as the data structure to store the matrices, leaving this aspect to the discussion of the practical implementation.
Thus, we make all worst-case complexity comparisons assuming the matrices are stored as dense columns.
While not perfect, this choice allows us to compare the pseudocodes in fairness, since the same disadvantage is applied to all.
When we compare the implemented versions, we will analogously choose for all implementations the best data structure (i.e., some sparse column representations), to ensure fairness.

\section{Background}
\label{sec_background}

In this section, we provide the formal definitions and constructions of the terms we will use.
The computation of PH is based on matrix operations, specifically multiplication and addition with matrices and columns.
The matrices we are interested in are constructed from the simplices (see \cref{subs_simplices}), and we restrict to the case in which they only have 0 and 1 as possible entries (that is, they have $\mathbb{Z}_2$ coefficients).
This allows us to avoid the notion of orientation in a simplex, and working modulo 2 makes addition and subtraction equivalent (that is, $1+1=0$, and $-1=1$).

\subsection{Matrix reduction}
Given a matrix $D$, we denote by $D[i]$ its $i$-th column, and by $D[i,j]$ the $i$th element of its $j$th column.
As in linear algebra, the \define{pivot} of $D[i]$, here denoted by $\piv{D[i]}$, is the row index of the lowest non-zero element in column $i$.
A matrix is \define{reduced} if all its non-zero columns have different pivots.
The standard PH \define{reduction} algorithm for making an unreduced matrix reduced is recalled in \cref{alg_sba}, here referred to as \nameref{alg_sba} - the Standard Barcode Algorithm \cite{standard_alg}.

\begin{figure}[h]
	\centering
	\input{figs/matrix-simplex.tex}
	\caption{A filtered simplicial complex $K$ (left) with its simplices indexed by entrance time, its unreduced, total boundary matrix (middle), and its barcode (right), constructed from the pivot indices of the reduced boundary matrix. In the barcode, dimension 0 bars are below the dashed line, and dimension 1 bars are above it.}
	\label{fig_matrix-simplex}
\end{figure} 

\subsection{Simplicial complexes and filtrations} 
\label{subs_simplices}
Given a finite set $V$, a \define{$k$-simplex} (over $V$) is a subset of $V$ of size $k + 1$.
Whenever there is a relationship $\tau \subsetneq \sigma$ among simplices, the simplex $\tau$ is called a \define{face} of $\sigma$ and $\sigma$ is called a \define{coface} of $\tau$.
If, in addition, $\dim(\tau)=\dim(\sigma)-1$, then $\tau$ is also called a \define{facet} of $\sigma$, and $\sigma$ is called a \define{cofacet} of $\tau$.
A \define{simplicial complex} is a set $K$ of simplices that is closed under non-empty subsets, that is, if $\tau \subset \sigma$ and $\sigma\in K$, then $\tau\in K$.
Throughout the paper, we assume the simplicial complexes are finite, that is, all simplices in the simplicial complex come from a finite set of vertices.
The \define{star} of a simplex $\sigma\in K$ is the collection of simplices containing $\sigma$, written $\str{\sigma}\coloneqq\{\tau\in K\ \colon\ \sigma \subseteq \tau \}$.
Note that $\str{\sigma}$ is in general not a simplicial complex, but $K\setminus \str{\sigma}$ is a simplicial complex for all $\sigma\in K$.
A simplicial complex $K$ has a naturally associated $n\times n$ \define{boundary matrix} $D$, whose $(i,j)$-entry $D[i,j]$ is $1$ if the $i$th simplex is a facet of the $j$th simplex, and $0$ otherwise, with simplices ordered by their entrance times.
Taking advantage of notation, we use the the same names for simplex indices  and column names in $D$.

A \define{filtration} of a simplicial complex $K$ (also called a \define{filtered simplicial complex}) is a sequence $\mathcal{F}\coloneqq K_1 \subseteq \cdots \subseteq K_n$ of simplicial complexes, with $K_n=K$.
Given an ordering $\preccurlyeq$ of all the simplices in a simplicial complex $K=\{\sigma_1,\dots,\sigma_n\}$ that respects the face relation (that is, if $\sigma_i\subseteq\sigma_j$ then $\sigma_i\preccurlyeq \sigma_j$), the \define{simplex-wise filtration of $K$} is defined by $K_{i}=\{\sigma_1,\dots,\sigma_i\}$.
This is the default filtration for our purposes.
The variable $n$ denotes the number of simplices in the simplicial complex, and thus the number of steps in the filtration.
The \define{entrance time} of simplex  $\sigma_i$, or step in the filtration when $\sigma_i$ first appears, is $i$.
The \define{(total) boundary matrix} of a filtration $\mathcal{F}$ is the boundary matrix of the simplicial complex at the last step of the filtration whose columns and rows are ordered by increasing entrance time.

\subsection{Homology}
We provide an intuition of homology and its computation, and refer the reader to \cite{dey_wang_2022,edels_harer,oudot} for formal details.

Given a simplicial complex $K$, we can consider sums of simplices of the same dimension, computed as sums of the respective columns in the boundary matrix.
For example, the sum $e_2+e_3$ of the simplicial complex in \cref{fig_matrixred}, indicated by dashed lines, is recorded in the boundary matrix by adding the second column to the third columns.
As a result, the \define{boundary} of the sum $e_2+e_3$ is $v_1+v_2$, as in \cref{fig_matrixred}(b), indicating that the boundary contains $v_1$ and $v_2$, which is exactly what we expect from the picture.
If we add also $e_1$, we obtain a closed path, which, intuitively, has no boundary and whose boundary column is indeed zero, as in \cref{fig_matrixred}(c).
A formal sum whose boundary column is zero is called a \define{cycle}, which also reflects the intuition of a closed loop, or ``hole''.

\begin{figure}[h]
	\centering
	\input{figs/matrix_reduction.tex}
	\caption{A simplicial complex (left), the sub-matrix of its boundary matrix given by the edges and their boundaries (a), the formal sums $e_2+e_3$ (dashed) and $e_1+e_2+e_3$ (dashed and solid), (b) and (c) respectively.}
	\label{fig_matrixred}
\end{figure}

If the simplicial complex in \cref{fig_matrixred} had also the $2$-simplex with edges $e_1,e_2,e_3$, then the formal sum $e_1+e_2+e_3$ would be its boundary.
Intuitively, if the $2$-simplex is there then there is no ``hole''.
This intuition gives us an informal definition of homology:
{\it (Nontrivial) \define{homology classes} are sums of simplices that are cycles, but not boundaries}.
Without such a 2-simplex, the cycle does not correspond to a boundary, and so is a representative cycle of a homological class, that is, a $1$-dimensional hole.
Note that a class may have several representative cycles, whenever their difference does not exist in homology, i.e., whenever they differ by a boundary.

\subsection{Persistent homology}
\label{sec_persistent_homology}
Given a filtration over a simplicial complex, its \define{persistent homology} is given by the homology of each complex in the filtration, together with morphisms connecting them.
When a simplex $\sigma$ enters the filtration, one of three things may happen:
\begin{enumerate}
	\item \label{case1} adding the column of $\sigma$ to an existing sum of simplices makes the boundary of the sum zero; or
	\item \label{case2} adding the column of $\sigma$ does not make any existing sum zero, and the boundary of $\sigma$ is not a cycle; or
	\item \label{case3} the boundary of $\sigma$ is the cycle of a formal sum of simplices entering earlier.
\end{enumerate}

\noindent\textbf{Case 1}: corresponds to the first occurrence in a filtration of a nontrivial homology class $c$ at the entrance time of $\sigma$.
This time is called the \define{birth time} of $c$ class, and $\sigma$ is called a \define{positive} simplex that is \define{born} at this entrance time and which \define{creates} the class $c$.

\noindent\textbf{Case 2}: there is no change in homology classes.

\noindent\textbf{Case 3}: corresponds to a filtration step without a homology class $c$ that existed at the previous step.
This time is called the \define{death time} of $c$, and $\sigma$ is called a \define{negative} simplex that is \define{killed} at this time and which \define{kills} the class $c$.
The interval between the birth and death times is called a \define{bar}, and the collection of bars is represented in a \define{barcode}, as in \cref{fig_matrix-simplex} (right).
In persistent homology, cycles that are not boundaries at some step (but may be boundaries later) are still called representative cycles.

The barcode of a filtration $\mathcal{F}$ is computed from its boundary matrix $D$, usually using (variations of) \nameref{alg_sba} \cite{standard_alg}.
The key part of any barcode computation is finding the \textit{pivot pairs}.
Let $D[\geq l, \leq f]$ denote the submatrix of $D$ given by the last $l$ rows of the first $f$ columns of $D$.
As shown in \cite{standard_alg}, $(i,j)$ is a \define{pivot pair} if and only if
\begin{align}
	\rank \left(D[\geq i,\leq j]\right) & - \rank \left(D[\geq i+1,\leq j]\right) \nonumber \\
	& + \rank \left(D[\geq i+1,\leq j-1]\right)   - \rank \left(D[\geq i,\leq j-1]\right) = 1. \label{eq_pivotpairs}
\end{align}

Given a pivot pair $(i,j)$ in a simplex-wise filtration, the simplex $\sigma_i$ creates a homological class and $\sigma_j$ kills that class, and we draw the interval from $i$ to $j$ as a bar in the barcode.
Most computations of the barcode \cite{adams2014javaplex, bauer2021ripser, phat_paper,henselman2016matroid, maria2014gudhi,dionysus, perez2021giotto} implement, with various speed-ups, the \nameref{alg_sba}, which we report here for self-containment with the slight modification of keeping explicit track (in \cref{line_opmat}) of the performed column operations in the matrix $V$.
As anticipated in \cref{subs_(pseudo)code}, we will not focus on how the columns are stored in a computer at this stage.\footnote{In the original paper, the \nameref{alg_sba} is presented using a \emph{linked list} data structure.}
We simply remark that, in practice, this algorithm is rarely applied as it is: the clearing optimization \cite{twist} (recalled in \cref{sec_relationtotwist}) is applied virtually always, often paired with cohomology \cite{ripser_software,dualities_persistent} to speed up the computation considerably.
In this work, we discuss the addition of the clearing (see \cref{sec_relationtotwist}), but we will not compare with the cohomology since we are interested in obtaining the representative cycles (and not the representative cocycles, which would be obtained using cohomology), leaving this aspect to future work.

\begin{algorithm2e}
	\caption[SBA]{\sc{SBA - \textbf{S}tandard \textbf{B}arcode \textbf{A}lgorithm} \cite{standard_alg}}
	\label{alg_sba}
	\KwIn{Boundary matrix $D$, identity matrix $I$}
	\KwOut{Reduced matrix $R$, operations matrix $V$}
	$R \leftarrow D$ \;
	$V \leftarrow \mathrm{I}$ \;
	\For{$j = 1,\dots,n$\nllabel{line_sbaforloop}}{
		\While{there exists $j'<j$ for which $\piv{\col{R}{j'}}=\piv{\col{R}{j}}\neq 0$}{
			add $\col{R}{j'}$ to $\col{R}{j}$ \;
			add $\col{V}{j'}$ to $\col{V}{j}$ \nllabel{line_opmat}}
	}
	\Return{$R$, $V$}
\end{algorithm2e}

The pivots of the resulting matrix $R$ are precisely the pivot pairs of the filtration given by $D$, hence the barcode can be read directly off $R$.
Note that \nameref{alg_sba} keeps track of the column operations used to reduce the matrix by storing them in $V$, giving the factorization $R=DV$.
To obtain all representative cycles, which are useful in practice as they describe the location of a hole in the data, one can use the matrix $V$, which also has the advantage of remembering all (column) operations performed during the reduction.
If, on the other hand, one is interested only in the representative cycles of the bars with finite length and does not need to recover the performed column operation, it is enough to record $R$.
Our method, presented in \cref{sec_removing,sec_necessity}, provides updates to both $V$ and $R$ when removing simplices from the filtration.
At the end of \nameref{alg_sba}, $V$ is upper triangular.
However, several variants of the SBA obtain a decomposition $R=DV$ where $V$ is not upper triangular, but $R$ still contains the correct pivot pairs \cite{keepitsparse}.
Even though other methods implementing barcode updates (such as \cite{vines_and_vineyards}) require $V$ to be upper triangular, we do not require it, mostly to avoid unnecessary extra computation steps.

\section{Removing simplices}
\label{sec_removing}

In this section, we enter in the core of our contribution: an efficient method to update a reduced boundary matrix when some of its simplices are removed.

Let $\mathcal{F}$ be a filtration of a simplicial complex $K$ and $\{\tau_1,\dots,\tau_t\}\subseteq K$ a set of simplices to be removed from $\mathcal{F}$.
Let $L=\bigcup_{i=1}^t \str{\tau_i}$ be the union of the stars of all the simplices $\tau_i$.
We now introduce a method to update the reduced boundary matrix of $\mathcal{F}$ to the reduced boundary matrix of the filtration $\mathcal{F}$ without the simplices in $L$, which we denote by $\mathcal F\setminus L$.
Note that, by definition of star, $\mathcal F \setminus L$ is a filtration, whereas $\mathcal F\setminus \{\tau_1,\dots,\tau_\ell\}$ may not be well-defined, as $K\setminus\{\tau_1,\dots,\tau_\ell\}$ may not be a simplicial complex.

Our main algorithm is \nameref{alg_rem} -- \textbf{Si}mplicial \textbf{R}emoval \textbf{U}pdate \textbf{P}rocedure (\cref{alg_rem}).
For the purposes of this algorithm, a zero column is considered to have $-1$ as a pivot, that is, the pivot of a zero column is higher than the one of any other non-zero column.

\begin{algorithm2e}
	\caption[\nameshort]{\sc{\nameshort - \textbf{Si}mplicial \textbf{R}emoval \textbf{U}pdate \textbf{P}rocedure}}
	\label{alg_rem}
	\KwIn{Reduced matrix $R'$, operations matrix $V'$, list of simplices to be removed $L$}
	\KwOut{Updated reduced matrix $R$, updated operations matrix $V$}
	$R \leftarrow R'$ \;
	$V \leftarrow V'$ \;
	\For{$\sigma_j \in L$\nllabel{line_for_simplicestoberemoved}}{
		$A = \{i\ :\ \col{V}{j,i}\neq 0,\ i\neq j\} = \{A_1,A_2,\dots\}$ in increasing order \nllabel{line_r} \;
		$B =\{A_0=j\}\cup \{A_i\in A\ :\ \piv{\col{R}{A_i}} < \piv{\col{R}{A_t}}\ \forall\ t < i\}$ in increasing order \nllabel{line_setb}\;
		\For{$i \in A$ in decreasing order \nllabel{line_for_affectedsimplices}}{
			$k\leftarrow$ the smallest index in $B$ such that $\piv{\col{R}{B_k}} \leqslant \piv{\col{R}{i}}$ \nllabel{line_analyze_pivot}\;
			\If{$i=B_k$ \nllabel{line_bcase}}{
				$k \leftarrow k-1$
			}
			add column $B_k$ to column $i$ in $R$ and in $V$ \nllabel{line_coladd}\;				
	}
	remove $\col{R}{j}$, $\col{V}{j}$, and row $j$ in $R$ and in $V$ \nllabel{line_removecolsrows}
}
\Return{$R$, $V$}
\end{algorithm2e}

The correctness of \nameref{alg_rem} is proved by \cref{thm_correctness}, using the correctness of \nameref{alg_rem_naive} (\cref{alg_rem_naive}).
The correctness of \nameref{alg_rem_naive}, that is, the claim that it produces the same pivot pairing as the \nameref{alg_sba} on $\mathcal{F}\setminus L$, is immediate, as \cref{line_for_simplicestoberemoved_naive} in \nameref{alg_rem_naive} undoes all column additions in which elements of $L$ were involved, and \cref{line_sba} ensures that the output is reduced.

\begin{algorithm2e}
\caption{\sc{Naive removal}}
\label{alg_rem_naive}
\KwIn{Reduced matrix $R'$, operations matrix $V'$, list of simplices to be removed $L$}
\KwOut{Updated reduced matrix $R$, updated operations matrix $V$}
$R \leftarrow R'$ \;
$V \leftarrow V'$ \;
\For{$\sigma_j \in L$\nllabel{line_for_simplicestoberemoved_naive}}{
	$A = \{i\ \colon\ \col{V}{j,i}\neq 0,\ i \neq j\}$ sorted in increasing order \nllabel{line_r_naive} \;
	\For{$i \in A$ in increasing order \nllabel{line_for_naive}}{
		$\col{V}{i} \leftarrow \col{V}{i} + \col{V}{j}$ \nllabel{line_addv_naive}\;
		$\col{R}{i} \leftarrow \col{R}{i} + \col{R}{j}$ \nllabel{line_addp_naive}\;
	}
	remove $\col{R}{j}$, $\col{V}{j}$, and row $j$ in $R$ and in $V'$ \nllabel{line_naive_remove} \;
	call \nameref{alg_sba}$(R,V)$\nllabel{line_sba}
}
\Return{$R$, $V$}
\end{algorithm2e}

The key difference between the two algorithms is that \nameref{alg_rem} skips repeated column additions, as made precise by \cref{lemma_bpivots,lemma_apivots}.
The addition of zero columns may also be skipped, when $\col{R}{B_k}$ is zero in \cref{line_coladd}, to reduce even more the number of operations.
The efficiency of \nameref{alg_rem} is proved in \cref{complexity_removal}, and demonstrated on practical examples in \cref{sec_experiments}.
An example to provide intuition on how \nameref{alg_rem} is built on top of \nameref{alg_rem_naive} is given by \cref{ex_sirup_worked_example}.

\begin{example}
\label{ex_sirup_worked_example}
We compare the steps of \nameref{alg_rem_naive} and \nameref{alg_rem} for the same removal.
In the presented matrices in this example, the darker squares (of any color) represent non-zero elements, and the lighter squares represent zero entries.
The red columns and rows correspond to the columns and rows of the simplex marked for removal, the yellow ones are the rows and columns being modified at a given step, and the blue parts are all the other entries.

Let $\mathcal{F}$ be the filtered simplicial complex depicted in \eqref{eq_sirup_worked_example_1}, with vertices in numerical order followed by edges in alphabetic order.
Let $R''=D''V''$ be the associated matrix reduction, as computed by \nameref{alg_sba}, with $L = \{c,a\} \subseteq K_1$.
The $R''=D''V''$ matrices restricted to relevant submatrices are given by
\begin{equation}
\input{figs/sirup_worked_example_new_1.tex}
\label{eq_sirup_worked_example_1}
\end{equation}

We remove $c$ first.
To do so, the rows and columns marked in red in \cref{eq_sirup_worked_example_1} will be removed.
Following \cref{line_r}, equally in \nameref{alg_rem} and \nameref{alg_rem_naive}, we get $A = \{d\}$, as row $c$ in $V''$ has only one non-zero off-diagonal entry.
\nameref{alg_rem_naive} asks to only add column $c$ to $d$, in both $R''$ and $V''$.
After removing $c$, the factorization is still reduced, so there is nothing for \nameref{alg_sba} to do.
Similarly, in \nameref{alg_rem} we have $B=\{c,d\}$, and as $9 = \piv{\col{R''}{c}} > \piv{\col{R''}{d}} = -1$, the only addition to do is also adding column $c$ to $d$.
We then get (in both algorithms) the updated matrices
\begin{equation}
\input{figs/sirup_worked_example_new_2.tex}
\label{eq_sirup_worked_example_2}
\end{equation}

To next remove $a$, the rows and columns marked in red in \cref{eq_sirup_worked_example_2} will be removed.
Construct $A=\{b,d,e,f,g,h,i,j,k\}$, as the row of $a$ in $V'$ has only non-zero elements.
For \nameref{alg_rem_naive}, after the for-loop in \cref{line_for_naive} and column removal in \cref{line_naive_remove}, the row of $V$ is empty off the diagonal, and $R$ is not reduced.
For \nameref{alg_rem}, the construction of $B$ on \cref{line_setb} identifies the earliest-highest pivots in $R$,
which are $\{a,b,f,i\}$ as emphasized in \cref{eq_sirup_worked_example_3}.
To keep track of the changes in the two different algorithms, we denote by $R_N=D_NV_N$ (resp. $R_S=D_SV_S$) the factorization as updated so far by \nameref{alg_rem_naive} (resp. by \nameref{alg_rem}), even if $V_S=V'$ and $R_S=R'$.
The updated matrices are given by
\begin{equation}
\input{figs/sirup_worked_example_new_3.tex}
\label{eq_sirup_worked_example_3}
\end{equation}

In \nameref{alg_rem_naive}, the \nameref{alg_sba} now proceeds to reduce the matrix $R_N$, from left to right.
The pattern of column additions to get a reduced matrix by \nameref{alg_sba} and its effect in $V_N$ indicates successive additions which cancel each other out.
\nameref{alg_rem} is based on this observation, effecting the single column addition which is the net result of the additions from \nameref{alg_rem_naive}.
The step-by-step operation can be seen in the left side of \cref{fig_sirup_worked_example_new_comparison}.

As shown in the right side of \cref{fig_sirup_worked_example_new_comparison}, in \nameref{alg_rem}, each column in $R_S$, from right to left, will have a column from $B$ added to it, determined by the relationship of the pivots: if the column is not in $B$, we find the earliest column in $B$ with pivot above the column we are considering.
If the column is in $B$, we take the column ordered just before the one we are considering.
This computation is completed in full detail in \cref{fig_sirup_worked_example_new_comparison}.
At the end, \nameref{alg_rem_naive} will have performed $9+12=21$ column additions, and \nameref{alg_rem} will only have performed 9.
\end{example}

\begin{figure}[htbp!]
	\centering
	\input{figs/sirup_worked_example_new_comparison.tex}
	\caption{Comparison from \cref{ex_sirup_worked_example} of the matrices $R$ (left side) and $V$ (right side) in \nameref{alg_rem_naive} (left column) and in \nameref{alg_rem} (right column).
	The label \textit{Column x} means that the algorithm is processing the $x$th column, and the arrows show column addition.
	Pivots in $B$ for \nameref{alg_rem} are emphasized.
	}
	\label{fig_sirup_worked_example_new_comparison}
\end{figure}

\subsection{Correctness and worst-case complexity}
\label{sec_correctness_complexity}
We now make rigorous the intuition from \cref{ex_sirup_worked_example} by proving the correctness of \nameref{alg_rem}.
We begin with two preliminary results about \nameref{alg_rem_naive} and the set $B$ defined in \cref{line_setb} of \nameref{alg_rem}.
For ease of notation, given a fixed $\sigma_j\in L$, we write $\col{R^\before}{i}$ for the state of $\col{R}{i}$ at the beginning of step $j$ in the for-loop of \cref{line_for_simplicestoberemoved_naive} of \nameref{alg_rem_naive}, and $\col{R^\after}{i}$ for the state of $\col{R}{i}$ at the end of that step, with analogous notation for $V$.

For convenience, the elements of $B$ from \cref{line_coladd} of \nameref{alg_rem} are indexed beginning with 0, that is, $B=\{B_0=j,B_1,B_2,\dots\}$.

\begin{lemma}
\label{lemma_bpivots}
In \nameref{alg_rem_naive}, for every $B_{i}\in B$ with $i>0$, we have
\begin{enumerate}
	\item $\piv{\col{R^\after}{B_i}} = \piv{\col{R^\before}{B_{i-1}}}$, and
	\item $\col{V^\after}{B_i} = \col{V^\before}{B_i} + \col{V^\before}{B_{i-1}}$.
\end{enumerate}
\end{lemma}

This lemma makes precise the cases in which \nameref{alg_rem_naive} changes a pivot. For example, in the extreme case that $B=A$, this lemma describes how
\begin{equation}
\input{figs/lemma_bpivots.tex}
\end{equation}

\begin{proof}
This follows by strong induction on the elements in $B$.

For the base case $i=1$, first note that $B_0=j$.
The column $\col{R}{B_1}$ is the left-most column with pivot higher than that of column $B_0$, and so it acquires the pivot of $\col{R}{B_0}$.
As $R$ was reduced, the pivot $\piv{\col{R}{B_0}}$ was unique among all previous pivots, so after removing $\col{R}{B_0}$, $\piv{\col{R^\after}{B_1}}$ is unique among all previous pivots.

For the inductive step, consider the column $\col{R}{B_i}$, to which we must add $\col{R}{j}$ in \cref{line_addp_naive} of \nameref{alg_rem_naive}.
In \cref{line_sba}, where SBA$(R,V)$ is executed, at the beginning of step $B_i$ of the for-loop in \cref{line_sbaforloop} of \nameref{alg_sba}, the inductive hypothesis on $R$ for $i=1$ gives us that column $\col{R}{B_1}$ is the only column with the same pivot as the column $\col{R}{B_i}$, and so $\col{R}{B_1}$ must be added to $\col{R}{B_i}$.
Now column $\col{R}{B_i}$ has the pivot $\piv{\col{R^\before}{B_1}}$, as the two earlier additions of $\col{R}{j}$ to each of $\col{R}{B_1}$ and $\col{R}{B_i}$ by \nameref{alg_sba} cancel each other out.
Again by the inductive hypothesis on $R$, now for $i=2$, $\col{R}{B_2}$ must be added to $\col{R}{B_i}$, as these two columns have the same pivot.
Now $\col{R}{B_i}$ has the pivot $\col{R^\before}{B_2}$, by the inductive hypothesis on $V$ for $i=2$.
This pattern continues until and including $\col{R}{B_{i-1}}$, as then $\col{R}{B_i}$ is the first column to have the pivot of $\col{R^\before}{B_{i-1}}$.
Then step $B_i$ of the for-loop in \nameref{alg_sba} finishes, as the pivot of $\col{R}{B_i}$ is now unique among all columns to its left.

By strong induction, the claim holds for all elements of $B$.
\end{proof}

Note that the second statement of \cref{lemma_bpivots} implies the first statement, but we make the distinction to emphasize the effect on pivots.
The second preliminary result follows a similar argument, while considering all other affected columns.

\begin{lemma}
\label{lemma_apivots}
In \nameref{alg_rem_naive}, for every $A_i\in A\setminus B$,
\begin{enumerate}
	\item $\piv{\col{R^\after}{A_i}} = \piv{\col{R^\before}{A_i}}$, and
	\item $\col{V^\after}{A_i} = \col{V^\before}{A_i} + \col{V^\before}{b}$,
\end{enumerate}
where $b$ is the first element in $B$ such that $\piv{R[b]}\geqslant \piv{\col{R^\before}{A_i}}$.
\end{lemma}

This lemma makes precise the cases in which \nameref{alg_rem_naive} does not change a pivot. For example, in the extreme case that $B=\{j\}$, this lemma describes how
\begin{equation}
\input{figs/lemma_apivots.tex}
\end{equation}

\begin{proof}
Consider some $A_i\in A\setminus B$, with $b=B_k$ as in the statement.

After $\col{R}{j}$ has been added by \cref{line_addp_naive}, \nameref{alg_sba} is executed in \cref{line_sba}.
We concentrate on the stage of \nameref{alg_sba} in which column $A_i$ is being considered, i.e. after all columns to the left of column $A_i$ have been reduced.
We claim that \nameref{alg_sba} performs precisely $k$ column additions to column $A_i$, with the column $\col{V^\after}{A_i}$ fully described by the telescoping sum
\begin{equation}
	\label{eqn_telescoping_v}
	\col{V^\before}{A_i} + \col{V}{B_0} + \col{V^\after}{B_1} + \cdots + \col{V^\after}{B_k} \equiv_2 \col{V^\before}{A_i} + \col{V^\before}{B_k} ,
\end{equation}
simplified using \cref{lemma_bpivots}.
This is straightforward if $k=0$, as then the addition of $\col{R}{j} = \col{R}{B_0}$ to column $A_i$ does not change the the pivot of column $A_i$.

For $k\geqslant 1$, the claim follows by observing that the pivot of $\col{R}{A_i}$ is now $\piv{\col{R}{B_0}}$, and so \nameref{alg_sba} will add $\col{R^\after}{B_1}$ to it.
Columns $\col{R^\after}{B_{k'}}$, for increasing $k'=1,2,\dots$, will continue to be added to $\col{R}{A_i}$, and the pivot of $\col{R}{A_i}$ will keep changing to $\piv{\col{R^\after}{B_{k'}}}$, until its pivot is unique (or it is the zero column).
This will occur after $k$ additions, as that is the first time that the sum in \cref{eqn_telescoping_v} leaves the column $\col{R}{A_i}$ with the pivot $\piv{\col{R^\before}{A_i}}$, as at this step $\piv{\col{R^\after}{B_{k+1}}}$ is finally higher (smaller index) than $\piv{\col{R^\before}{A_i}}$.
The pivot $\piv{\col{R^\before}{A_i}}$ is unique, as the input matrix was reduced.

Hence the claim for \nameref{alg_sba} holds, and the consideration of column $A_i$ ends after $k$ additions.
\end{proof}

These two results allow us to prove our main statement.

\begin{theorem}
\label{thm_correctness} 
Let $\mathcal{F}$ be a filtration of a simplicial complex $K$ and $L=\bigcup_{i=1}^t \str{\tau_i}$, where $\{\tau_1,\dots,\tau_t\}\subseteq K$ is a set of simplices to be removed from $\mathcal{F}$.
Using a reduced matrix and operations matrix for $\mathcal F$ as input, the output of \nameref{alg_rem} is a reduced matrix and operations matrix for $\mathcal F\setminus L$.
\end{theorem}

\begin{proof}
First note that the only modification of the input reduced matrix in \nameref{alg_rem} is done in \cref{line_coladd}, which is inside a for-loop iterating over all $A_i$ in decreasing order.
Fix $A_i\in A$.
If $A_i\in B$, by \cref{lemma_bpivots} the modification is the same as in \nameref{alg_rem_naive}.
If $A_i\in A\setminus B$, by \cref{lemma_apivots} the modification is the same as in \nameref{alg_rem_naive}.
Indeed, as the for-loop in \cref{line_for_affectedsimplices} is done in reverse, the column that is added to $\col{R}{A_i}$ is $\col{R^\before}{b}$, not $\col{R^\after}{b}$, as the index of $b$ in $A$ is necessarily smaller than the index $i$ of $A_i$.
This follows by the definition of $b$ from \cref{line_analyze_pivot} and the special case for elements of $B$ in \cref{line_bcase}.

Hence the modifications of the reduced boundary and operations matrices by \nameref{alg_rem} precisely match those of \nameref{alg_rem_naive}, and as the correctness of \nameref{alg_rem_naive} has been established, the statement holds.
\end{proof}

Note that the above arguments can be be extended to fields over a prime $p\neq 2$ (i.e., to matrices whose entries are not necessarily just 0 and 1), replacing column additions with column subtractions.
Indeed, the skipping of repeated operations hinges on whether an element is zero or not, rather than on its exact value.

\begin{proposition}
\label{complexity_removal}
Let $n$ be the number of simplices in the filtration $\mathcal F$ and $m$ the number of simplices in $L$.
Then \nameref{alg_rem} has worst case complexity $O(mn^2)$.
\end{proposition}

\begin{proof}
The for-loop in \cref{line_for_simplicestoberemoved} is executed $m$ times.
Constructing the list $A$ in \cref{line_r} requires at most $n$ operations, as each column of $V$ (and $R$) has to be considered once.
The construction of list $B$ in \cref{line_setb} may be done at the same time, as the pivot of every column has to be compared with the highest (smallest index) pivot seen so far in the construction of $A$.
If it is higher (smaller index), then that column index is added to $B$.
Identifying $k$, and so identifying $b$, in \cref{line_analyze_pivot} has complexity $O(|B|)=O(n)$, as at most every element in $B$ has to be considered.
Adding each column in \cref{line_coladd} is an $O(n)$ operation.
The removal operations in \cref{line_removecolsrows} require a constant number of operations, independent of $m$ or $n$.
Hence the worst case complexity is $O(mn^2)$.		
\end{proof}

Conversely, the complexity of computing the barcode of $\mathcal F\setminus L$ with \nameref{alg_sba} is $O((n-m)^3)$.
Note that this is the worst-case complexity, and does not take into account the choice of data structure.
In \cref{sec_experiments} we compare the experimental running times of \nameref{alg_rem} and implementations of \nameref{alg_sba}\cite{phat_paper}, demonstrating that the former is always faster when recomputing from scratch and computing the representative cycles, and almost always faster if one is only interested in computing the barcode.

Note that the size of $A$ is bounded above by the number of simplices in dimension $\dim(\sigma_j)$.
In particular, it is usually the case that $|A|\ll n$, that is, each column is added only to a small number of other columns.
We also note that if $|A|=0$, then $\sigma_j$ must be a maximal simplex (that is, a simplex without cofaces).

\subsection{Minimality}
In addition to being a valid algorithm that updates the boundary and operations matrices, \nameshort~provably executes the least amount of column additions necessary to keep the representative cycles of a filtration updated after a simplex $\sigma$ has been removed.
This follows as exactly one column addition is done for every cycle in which $\sigma$ appears, or equivalently, for every column in the operations matrix in which the row of $\sigma$ is nonzero.

It may be true that some of these column additions may not change the barcode (as is the case for column $j$ in \cref{ex_sirup_worked_example}), but without these column additions, the cycle representatives for non-trivial homology classes in the updated filtration can not be identified from the operations matrix.
Even more, not performing these column additions may lead to incorrect identification of affected simplices if the process is repeated on the output.

\subsection{Relation to clearing} 
\label{sec_relationtotwist}
A standard optimization of \nameref{alg_sba} called \define{clearing} \cite{twist, clearcompress} replaces columns of $R$ with zero columns, instead of performing the column operations that would ultimately reduce them to zero.
This optimization, based on the fact that the boundary of a boundary is zero, is commonly used in barcode algorithms, as it speeds up the computations considerably, especially when combined with cohomology \cite{twist,phat_paper,bauer2021ripser} or by doing row additions \cite{notes_pivot}.
An example of the effect is described in \cref{fig_clearing}.

\begin{figure}[h!]
	\centering
	\input{figs/clearing.tex}
	\caption{
	A filtered simplicial complex $K$ and its unreduced boundary matrix $D$ (left), factored as $R=DV$ by \nameref{alg_sba} (right, top), and the outputs $R_c,V_c$ with clearing (right, bottom).
	Clearing puts column $i=5$ of the boundary matrix $R_c$ to zero if and only if there is a pivot pairing $(i,j)=(5,6)$.
	That is, if the formal sum corresponding to column $j=6$ kills the homological class generated by simplex corresponding to row $i=5$.
	}
	\label{fig_clearing}
\end{figure}

This optimization is compatible with our method and may be selected in our implementation of \nameref{alg_rem} \cite{phat_vineyards} with the flag \texttt{--twist}.
Note that by default, the operations matrix $V$ would not encode the column operations that would reduce column $i$ to zero if that column is replaced by zero when clearing is used.
To still encode this information while using the clearing speedup, our implementation takes as a representative cycle of column $i$ the simplices in the boundary of column $j$, where $\col{R}{i,j}$ is the pivot triggering the clearing operation.
The described approach may be used to produce an updated factorization $R=DV$ in other dynamic settings in which the content of $V$ is pertinent, for example, when swapping the order of simplices \cite{vines_and_vineyards}.
That is, the clearing algorithm does not guarantee an updated factorization, and only updates $R$.

\section{Understanding removal}
\label{sec_understanding}
We now analyze different aspects of removing simplices at the geometric and barcode level, discuss the necessity and implications of updating $V$, and compare our method to the ones in the literature.

\subsection{Updating the barcode}
\label{sec_updatebarcode}
When two simplices swap order in a filtration, there are well-understood changes that happen to the barcode, fully described in \cite{vines_and_vineyards}.
Conversely, when a simplex and its cofaces are removed, it is not as straightforward to succinctly describe what happens to the barcode, as each of the cofaces is either creating or destroying a homological class.
Note that this is also the case when we insert a simplex and its cofaces.
Some examples of simplex removal and the effects on the barcode are demonstrated in \cref{fig_bigchanges}.

\begin{figure}[h]
	\centering
	\input{figs/bigchanges.tex}
	\caption{
	Examples of how the barcode may change when the star of a simplex (highlighted in red) is removed from different 2-dimensional simplicial complexes.
	Finite bars (one in dimension $0$ and several in dimension $1$) disappear and an infinite bar in dimension $1$ appears in (a), a finite bar becomes shorter and another finite bar disappears in (b), $4$ infinite bars in dimension $1$ appear in (c), and an infinite bar disappears in (d).
	Filtrations respect dimension, with all 0-simplices appearing before any 1-simplex, and all 1-simplices appearing before any 2-simplex.
	}
	\label{fig_bigchanges}
\end{figure}

Note that \cref{fig_bigchanges}.(c) may be generalized from 5 to $N$ 2-simplices as cofaces, creating $N-1$ infinite bars in dimension $1$.
This construction can realized in higher dimensions as well, creating ``fans'' of $k$-simplices whose maximal face is removed when the central facet, a shared $(k-1)$-simplex, is removed, triggering a change that adds $N-1$ infinite bars in dimension $k-1$.
If $\tau$ that is removed gave birth to a homological class, the dimension of that class must be $\dim(\tau)$, and if it kills a class, it can only kill a class of dimension $\dim(\tau)-1$.

\begin{remark}\label{claim_barcodechanges}
Let $\mathcal F$ be a filtration of the factorization $R=DV$ of a simplicial complex containing a simplex $\tau$, and $\mathcal F'$ the same filtration with $\tau$ removed.
Exactly one of the following must occur:
\begin{enumerate}
	\item \label{barcode_case1} The simplex $\tau$ is positive (its reduced column in $R$ is a zero column). Then there is an infinite bar in the barcode of $\mathcal F$ that disappears in the barcode of $\mathcal F'$.
	\item \label{barcode_case2} The simplex $\tau$ is negative and is not part of the representative cycle of any non-trivial homology class. Then the barcode of $\mathcal F$ has a finite bar killed by $\tau$, which becomes an infinite bar in the barcode of $\mathcal F'$.
	\item \label{barcode_case3} The simplex $\tau$ is negative and is part of the representative cycle of some  non-trivial homology class. Then the barcode of $\mathcal F$ has a finite bar killed by $\tau$, which is longer in the barcode of $\mathcal F'$ and is now killed by a coface of $\tau$. In addition, the bar in the barcode of $\mathcal F$ corresponding to this coface disappears in the barcode of $\mathcal F'$.
\end{enumerate}

Recall that the representative cycle of a (non-trivial) homology class consists of the simplices associated to the non-empty rows of the column in $V$ at which that class is born.
\end{remark}

The cases in \cref{claim_barcodechanges} can be tracked in \nameref{alg_rem}.
In case \ref{barcode_case1}, the pivot of the column of the simplex $\tau$ is higher (smaller index) than all other pivots, so $B$ from \cref{line_setb} is the singleton $\{j\}$, and $k=0$ for every step in the for-loop of \cref{line_for_affectedsimplices}, and the if-condition in \cref{line_bcase} is never satisfied.
In case \ref{barcode_case2}, the set $A$ of affected simplices from \cref{line_r} is the singleton $\{j\}$, as is $B$, so the for-loop of \cref{line_for_affectedsimplices} does not have any steps.
In case \ref{barcode_case3}, the set $A$ has at least two elements, the set $B$ has at least one element, and the for-loop of \cref{line_for_affectedsimplices} has at least one step.

\subsection{The necessity of \texorpdfstring{$V$}{k}}
\label{sec_necessity}
The \nameref{alg_rem} algorithm requires keeping track of class representatives in the matrix $V$, which is more than the minimal amount of information necessary to compute the barcode.
This is an advantage as it keeps the representative cycles updated, but it is more memory-consuming than just updating the persistence diagram.
One may ask if it is necessary to remember $V$ for barcode updates, in cases when representative cycles are not needed.
The answer is yes, as \cref{fig_reduced-matrices} demonstrates.
In this example it is not possible to update the barcode after a simplex has been removed by only remembering $R$, and $V$ contains the information that makes the update possible.
In \cref{fig_reduced-matrices}, the resulting matrix $R$ is the same, but different 1-simplices in different filtrations contributed to zeroing out the same column of $e_4$.
That is, the effect on the barcode computed by $R$ by removing a simplex is ambiguous, and by checking the matrix $V$, we can obtain the correct updated barcode.

We note that the representative cycle of a finite bar may be computed only with $R$, and that there are special cases in which the removal of a simplex only affects finite bars (see \cref{fig_bigchanges}).
As we see here and throughout, the matrix $V$ is necessary in general to update both the barcode and representative cycles after a simplex removal.
\begin{figure}[h]
	\centering
	\input{figs/reduced-matrices.tex}
	\caption{A simplicial complex without the 1-simplex $e_4$ (left) and three possible matrices $U_i$ for three different choices of $e_4$. For $D_i$ the appropriate boundary matrix in each case, all three choices satisfy $\left[\begin{smallmatrix} 0 & R \\ 0 & 0 \end{smallmatrix}\right] = D_i \left[\begin{smallmatrix} I & 0 \\ 0 & U_i \end{smallmatrix}\right]$.}
	\label{fig_reduced-matrices}
\end{figure}

\subsection{Comparison with other methods}
\label{sec_comparison}

\subsubsection{Removal as a dynamic zigzag}
\label{sec_zigzagcomparison}

Removing a simplex from a filtration can also be realized by expanding the filtration by adding backward maps at the end, thus turning it into a zigzag, where each backward map removes one simplex.
Dey and Hou \cite{dey2022updating,dey2023computing,fast_zigzag} develop an efficient method (called \deyhoualgo) for computing the persistent homology of zigzag filtrations, defined for a larger class of changes than just simplex removal.
Their method is very efficient compared to other zigzag computations, but, unsurprisingly since it is constructed for a more general setting, it is less efficient in practice (even if not asymptotically) than \nameref{alg_rem} for removing simplices.
Indeed, as detailed in \cref{ex_dey-hou-comparison}, for each simplex that needs to be removed, \deyhoualgo~performs $N(2N-1)$ operations.
On the other hand, \nameref{alg_rem} needs at most $N^2$ operations in general, and $0$ in the described example.

\begin{example}
\label{ex_dey-hou-comparison}
Consider a simplicial complex $K$ on $N$ vertices connected by $N-1$ edges, arranged as a line graph.
Dey--Hou requires simplex-wise filtrations which start and end at $\emptyset$, but for brevity we consider the truncation of only forward arrows.
The omitted right side in \cref{eq_deyhouexample0} is the mirror image of the left side (that is, backward arrows in reverse order of how the simplices were introduced).
An example filtration for $N=4$ is given by
\begin{equation}
	\label{eq_deyhouexample0}
	\input{figs/diag_dey-hou-comparison0.tex}
\end{equation}

Removing the edge $e_1$ is a backwards map $\xhookleftarrow{e_1}$.
This map already appears in the truncated end of \cref{eq_deyhouexample0}, and as $e_1$ is maximal in $K$, we can move it via backward and outward switches to immediately after the addition of $e_1$.
This is done so that when computing the barcode, ignoring all bars with endpoints at the introduction or removal step of $e_1$ will produce the barcode as if $e_1$ never existed.
The associated filtration $\hat {\mathcal E}$ required by Dey--Hou of only forward arrows is constructed by replacing a backwards arrow $\xhookleftarrow{\sigma}$ with a forwards arrow $\xhookrightarrow{\widehat\sigma}$ of the cone of the same simplex, coned at the common cone point $\omega$ (which has been added before any other simplex), in reverse order of all the backward arrows.
Note that the backwards arrow $\xhookleftarrow{e_1}$ must be moved again to, to have an ``up-down" filtrations.
Beginning at the step at which all the simplices of $K$ have been added, this new filtration of a simplicial complex $\widehat K$ ends in
\begin{equation}
	\label{eq_deyhouexample1}
	\input{figs/diag_dey-hou-comparison1.tex}
\end{equation}

At this step the filtration of only forward arrows may be reduced directly.
The filtrations of the simplicial complexes $K,\widehat{K}$ have reduced factorizations $R=DV$ and $\widehat R=\widehat D\widehat V$, respectively, where
\begin{equation}
	\label{eq_deyhouexample2}
	\input{figs/diag_dey-hou-comparison2.tex}\ ,
\end{equation}
with colors of emphasized simplices $e_1,\widehat{e_1}$ consistent with \cref{eq_deyhouexample0,eq_deyhouexample1}.
\nameref{alg_rem} requires 0 column additions to remove $e_1$ from the decomposition $R=DV$, as the column of $e_1$ is not added to any other column, as evidenced by $V$.
Computing $\widehat{R}$ by the standard barcode algorithm requires $N(2N-1)-1$ more column additions than computing $R$, as evidenced by $\widehat{V}$.
\end{example}

\subsubsection{Removal as a sequence of vine swaps}
\label{sec_vinescomparison}

The method of \cite{vines_and_vineyards} is adapted by \cite{luo2023accelerating} to perform a sequence of transpositions in the order of the simplices, moving a simplex that needs to be removed to the end of the filtration and then dropping it.
Compared with this approach, \nameref{alg_rem} does not require the operations matrix $V$ to be upper triangular, but we do require updates of $V$ to encode the correct representative cycles.
Contrary to what happens with the zigzags case in the previous section, swapping the vines will keep track of the representative cycles.
However, it is not as efficient as \nameref{alg_rem}.
In general, \cite{vines_and_vineyards} requires a single column (and row) addition in $R$ (and $V$) for every simplex with $\col{V}{\sigma,i}=1$, given by Cases 2.1 and 3.1 in \cite{vines_and_vineyards}.
This corresponds to every element (excluding $\sigma$) in the list $A$ in \nameref{alg_rem}, and similarly a single addition (of not necessarily the same summands) in each of $R$ and $V$ are performed in each step of the for-loop in \cref{line_for_affectedsimplices} of \nameref{alg_rem}.
Therefore, for every column addition performed by \nameref{alg_rem}, swapping vines has at least one column operation, and thus cannot be more efficient.
Moreover, there is at most one more column and row addition required by \cite{vines_and_vineyards} for each swap, for when pivots change (Cases 1.1.2, 2.1.2 in \cite{vines_and_vineyards}) as a result of the swap.
In contrast, no more additions are ever required by \nameref{alg_rem}.
Hence the method of \cite{vines_and_vineyards,luo2023accelerating} takes at least as many steps as \nameref{alg_rem}, to produce the same updated persistence pairs, with a potentially different $R=DV$ factorization.

In practice, swapping vineyards may incur $O(n)$ more operations than \nameref{alg_rem}.
For a practical example, consider the filtration in \cref{fig_bigchanges}(c), with $N=5$ simplices in dimension 2.
To remove the star of the 1-simplex, \nameref{alg_rem} performs one column addition in $R$ and $N$ column additions in $V$, as all of the columns to which the 1-simplex was added are zero.
Conversely, \cite{vines_and_vineyards} would perform $N+1$ column additions in both $R$ and $V$, as well as $2N$ column transpositions in $R,V$ for the same task.

To further illustrate the difference, we can use the filtration described in \cref{ex_dey-hou-comparison}.
We can use a sequence of swaps to move $e_1$ and $\widehat{e_1}$ to the end of the filtration $\widehat{\mathcal E}$.
This would require $2(N-1)$ column additions, as the algorithm \nameref{alg_sba} adds the column of $e_1$ to each of the columns of the coned vertices exactly once (except the first one), corresponding to Case 3.1 in \cite{vines_and_vineyards} and producing two column operations each.

\section{Implementation and testing}
\label{sec_experiments}

\begin{figure}[h]
\centering
\setlength{\tabcolsep}{10pt}
\begin{tabular}[t]{ccc}
	\includegraphics[scale=.6]{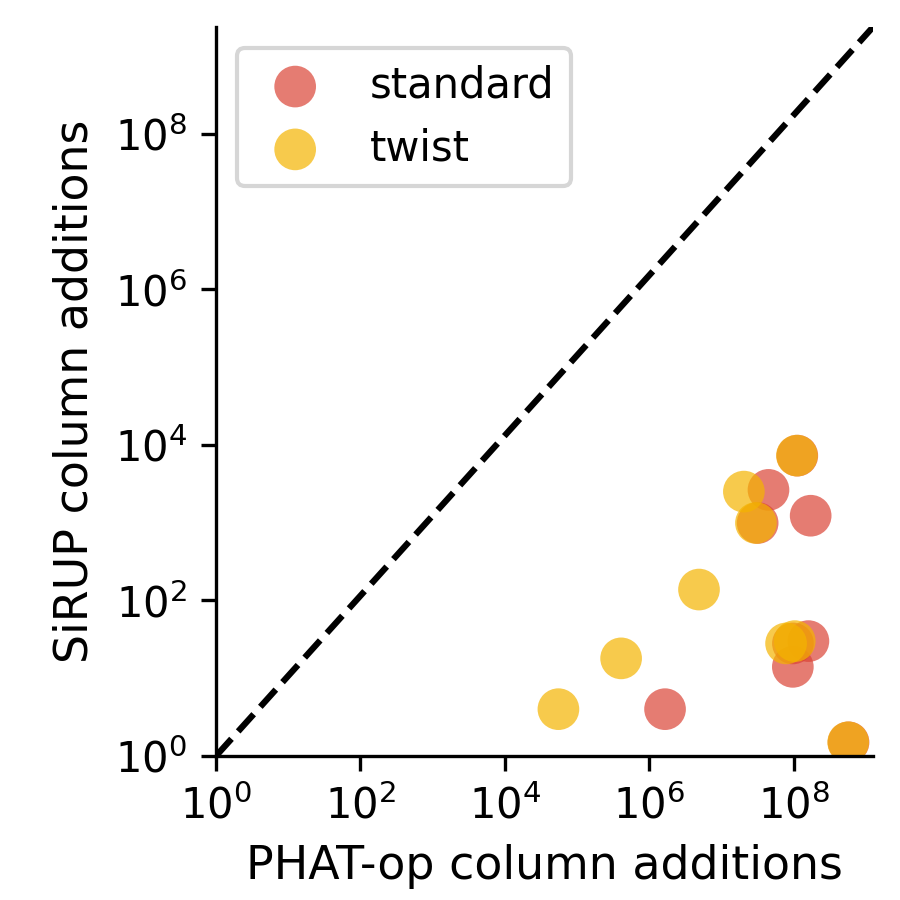} &
	\includegraphics[scale=.6]{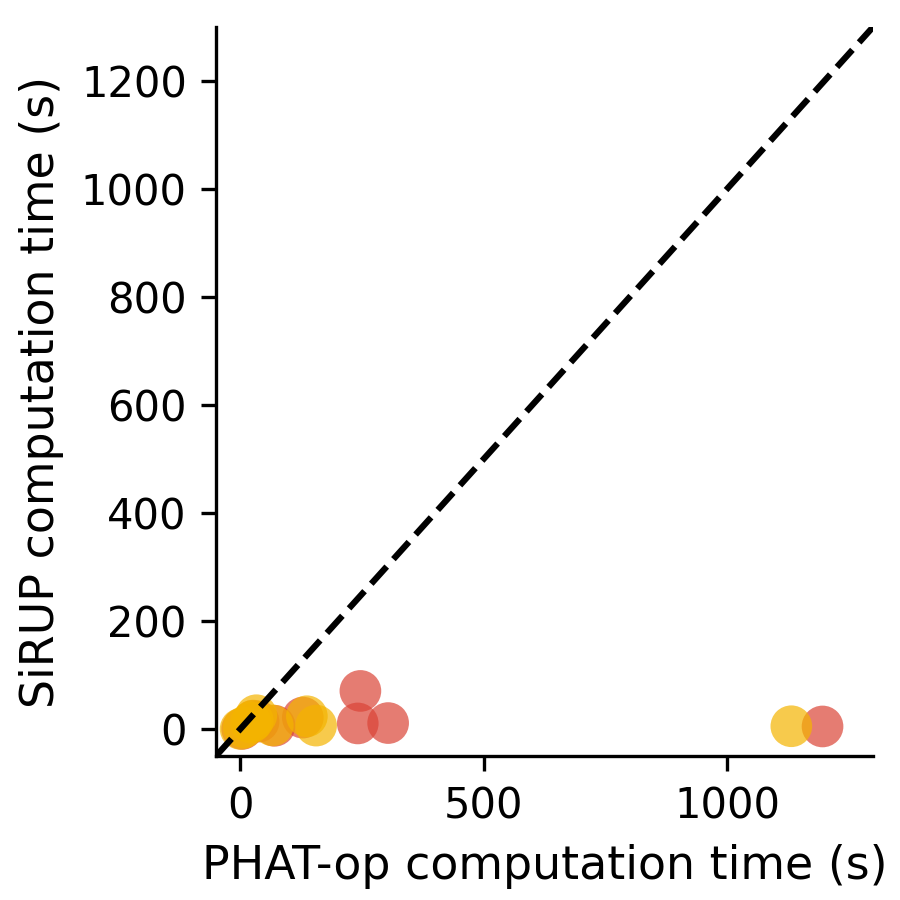} &
	\includegraphics[scale=.6]{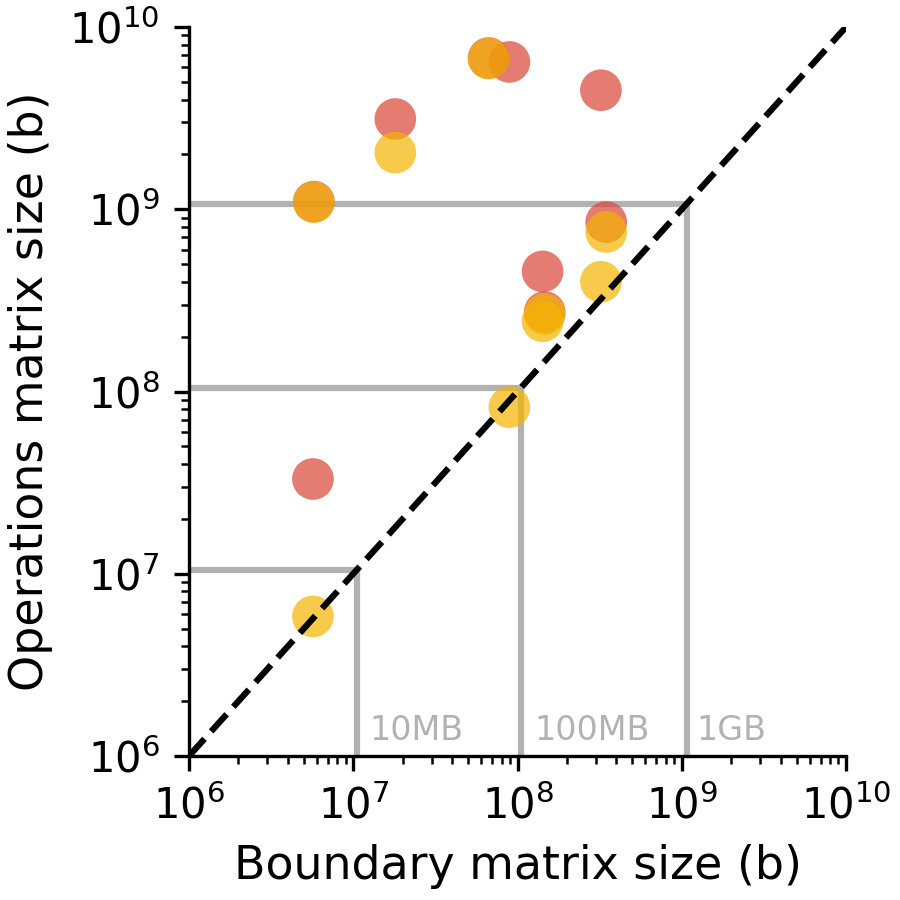}
\end{tabular}
\setlength{\tabcolsep}{6pt}
\caption{A visual overview of our main testing results.
	The plots are colored semi-transparently to distinguish the cases when there is little difference between the \texttt{--standard} and \texttt{--twist} methods. File sizes of boundary and operations matrices are also compared (right), highlighting the tradeoff in time for space for \nameshort. The full details of the left and center plots are given in \cref{table_updates}. Raw data for all three plots is given in the related dataset \cite{our_zenodo}.}
\label{fig_results}
\end{figure}

We have implemented \nameref{alg_rem} in a publicly available repository \phatv \cite{phat_vineyards}, and have performed tests to demonstrate its feasibility (\cref{sec_testing}).
We performed several experiments on a diverse set of inputs, recording both the time necessary to complete the experiments and the number of operations needed, with code to generate and execute the tests described also publicly available on Zenodo \cite{our_zenodo}.
Our main results, presented visually in \cref{fig_results} and in more detail in \cref{table_updates}, demonstrate that updating representative cycles and the barcode with \nameshort~is always faster and takes less matrix operations than doing the same with PHAT-op (PHAT software that also keeps track of the operations matrix).
Even when comparing with a task that does less, only computing the barcode from scratch with PHAT, \nameshort~is still almost always faster on tests that remove up to 200 simplices.

\subsection{Implementation}
\label{sec_implementations}

The implementation of \nameshort, named \phatv, is built on top of the methods of the Persistent Homology Algorithm Toolbox (PHAT) \cite{phat_paper}, a collection of \texttt{C++} headers to compute persistence pairs of a simplicial complex filtration.
This software was chosen because we wanted to compare just the reduction routine, which is precisely the focus of PHAT.
Following the findings of \cite{phat_paper} of efficiency in practice, we chose the \define{bit tree} as a data structure to store and perform operations with the columns of the matrices.
Nevertheless, any software that takes in a list of simplices and computes persistence pairs by way of operations on the boundary matrix is amenable to an implementation of \nameshort.
While the current implementation of \phatv~is already more efficient in most cases than recomputing from scratch, it is still a proof-of-concept: the look-up operations in \phatv~have potential to be be considerably optimized with a more suitable data structure.
However, our method is not yet adapted to handle cohomology.
Thus, if compared with PHAT when using the twist and cohomology optimizations, \nameref{alg_rem} would most likely perform worse, and would need to be adapted to be competitive in this case too.
We stress that, in all the implementations, we took advantage of the data structures implemented in PHAT to maximize efficiency.
In particular, none of the implementations (neither ours, nor the methods we compare with) uses dense columns, as this choice is much slower in practice.
Moreover, a proper implementation should include not just removal but integrate the swaps and additions of simplices from the literature, which goes beyond the scope of the current paper.
Here, the focus has been to accurately reflect \nameref{alg_rem} (\cref{alg_rem}) and demonstrate its role as a theoretical improvement over the current state-of-the-art.

\subsection{Testing}
\label{sec_testing}

Our experiments have been run on a computer with an Intel Core i7-8650U CPU with eight 1.9GHz cores and 15.4GiB of RAM, running Ubuntu 22.04 LTS, with gcc version 11.4.0.
The implementation is not parallelized.
Random selection of the simplices was done with \texttt{std::random\_device} seeding the pseudo-random generator \texttt{std::mt19337}.
Timing of the individual parts of the code was done with \texttt{omp\_get\_wtime}, and of the full execution by reading the \texttt{SECONDS} variable within the \texttt{bash} command interpreter of Linux.

\subsubsection{Inputs}
The datasets are either taken from \cite{roadmap,reliability,tadasets,images_database} or generated randomly with \cite{numpy}.
They are chosen to be as large as possible, while still being computationally manageable, and to demonstrate the strengths and weaknesses of \nameref{alg_rem}.
A full overview of the datasets and associated topological constructions is given in \cref{table_inputs,table_filtrations}.
\begin{table}[h]
\centering
\small
\begin{tabular}[t]{l|rrrrr|c}
	Dataset & 0-$\dim$ & 1-$\dim$ & 2-$\dim$ & 3-$\dim$ & 4-$\dim$ & Source \\ \hline
	\texttt{vr\_dsphere} & 3\ 000 & 65\ 052 & 400\ 187 & 1\ 224\ 795 & 2\ 353\ 767 &
	\cite{tadasets,flagser} \\
	\texttt{vr\_torus} & 1\ 000 & 21\ 988 & 214\ 105 & 1\ 415\ 261 & 7\ 288\ 105 &
	\cite{tadasets,flagser} \\
	\texttt{vr\_senate} & 103 & 5\ 253 & 176\ 851 & 4\ 421\ 275 & 0 &
	\cite{roadmap,flagser} \\
	\texttt{er\_clique} & 100 & 4\ 463 & 118\ 583 & 2\ 111\ 092 & 0 &
	\cite{numpy,flagser} \\
	\texttt{er\_shuffled} & 100 & 3\ 005 & 36\ 126 & 194\ 908 & 0 &
	\cite{numpy,flagser} \\
	\texttt{alpha\_cube} & 10\ 000 & 76\ 644 & 133\ 170 & 66\ 525 & 0 &
	\cite{tadasets,keepitsparse} \\
	\texttt{alpha\_swissroll} & 100\ 000 & 1\ 043\ 027 & 1\ 829\ 784 & 886\ 756 & 0 &
	\cite{tadasets,keepitsparse} \\
	\texttt{cc\_tooth} & 1\ 558\ 802 & 4\ 635\ 007 & 4\ 593\ 966 & 1\ 517\ 760 & 0 &
	\cite{images_database,keepitsparse}\\
	\texttt{bio\_bbmcl6} & 3\ 174 & 180\ 032 & 529\ 470 & 157\ 900 & 0 &
	\cite{reliability,flagser} \\
	\bottomrule
\end{tabular}
\caption{Descriptive overview of input datasets, with the number of simplices by dimension, and sources for sampling and construction indicated. The filtrations on each dataset that result in the complexes with the given numbers of simplices are described in \cref{table_filtrations}.}
\label{table_inputs}
\end{table}

We used five different types of filtrations (\define{Vietoris--Rips}, \define{alpha}, \define{random}, \define{lowerstar}, and \define{biological}) on two different types of complexes (\define{simplicial} and \define{cubical}), chosen for their common use in practice.
We will not define these filtrations here, but the unfamiliar reader may look the definitions up in \cite{roadmap} (the first four) and \cite{bb} (the latter).
The \define{clique} filtration was used together with the random filtration, by first randomly selecting the order of edges, then including simplices for all resulting cliques, also in a random order.
All complexes were simplicial complexes, except for \texttt{cc\_tooth}, which was a cubical complex defined by a three-dimensional scan of a tooth.
Random sampling for simplicial complexes based on Erd\H{o}s--R\'enyi graphs is done with NumPy \cite{numpy}.
Identifying higher dimensional simplices for Vietors--Rips and and Erd\H{o}s--R\'enyi constructions is done with the \textsc{Flagser-Count} \cite{flagsercount} variant of \textsc{Flagser} \cite{flagser}.
Full descriptions of the filtrations are given in \cref{table_filtrations}.

\begin{table}[h]
\centering
\footnotesize
\begin{tabular}[t]{p{.17\textwidth}|p{.10\textwidth}|p{.64\textwidth}}
	Dataset & Filtration & Description \\ \hline
	\texttt{vr\_dsphere} & Vietoris--Rips & \multirow{7}{*}{\parbox{.64\textwidth}{\raggedright
			Random samples over a unit 4-sphere in $\mathbb{R}^6$, random samples over a 2-torus with major radius 2 and minor radius 1 in $\mathbb{R}^4$, and a legislator voting dataset (distances among 103 points) from \cite{roadmap}. Both random samples have normally distributed noise at standard deviation $0.2$. Only edges less than 0.6 are considered for \texttt{vr\_dsphere}, and less than 1 for \texttt{vr\_torus}. All edges are considered for \texttt{vr\_senate}. }} \\
	\texttt{vr\_torus} & & \\ & & \\
	\texttt{vr\_senate} & & \\ & & \\ & & \\ \hline
	\texttt{er\_clique} & random, clique & \multirow{7}{*}{\parbox{.64\textwidth}{\raggedright
			Random values in $[0,1]$ assigned to $\binom{100}{2}$ edges, included from highest to lowest ``probability," with higher dimensional simplices defined by the clique complex. For \texttt{er\_clique}, a $(d>1)$-simplex is ordered immediately after all its faces. For \texttt{er\_shuffled}, the order of every $(d>1)$-simplex is included in a random order after all its faces have been included. Probability at least $0.1$ considered for \texttt{er\_clique} and at least $0.5$ for \texttt{er\_shuffled}. }} \\
	\texttt{er\_shuffled} & & \\ & & \\ & & \\ & & \\ & & \\ & & \\ \hline
	\texttt{alpha\_cube} & alpha & \multirow{2}{*}{\parbox{.64\textwidth}{\raggedright
			Random samples of 10\ 000 and 100\ 000 points over a cube and swiss, respectively, generated with \cite{tadasets}. }} \\
	\texttt{alpha\_swissroll} & & \\ \hline
	\texttt{cc\_tooth} & lowerstar, clique &
	A 3D image of a tooth from \cite{images_database}, represented as a $104\times 91 \times 161$ array of values in $[0,255]$. This defines a filtration of a cubical complex in 3 dimensions. \\ \hline
	\texttt{bio\_bbmcl6} & biological, clique &
	Biologically motivated reconstruction of a rat's neocortex \cite{bb}, with neurons as vertices, synapses as edges, and their clique complex defining higher dimensional simplices. Restricted to the undirected excitatory subcircuit of bipolar pyramidal cells in Layer 6 ($\approx 10\%$ of all the neurons). \\
	\bottomrule
\end{tabular}
\caption{Details of filtrations on input datasets.}
\label{table_filtrations}
\end{table}

\subsubsection{Updates}
For each dataset, a fixed number of simplices was selected to be removed in one of the dimensions of the complex.
The union of the stars of those simplices was computed, producing an update file with the indices of all the simplices to be removed, and an updated boundary matrix without those simplices.
The number of simplices, by dimension and in total, is reported in the second column of \cref{table_updates}.
\begin{table}[h]
	\centering
	\scriptsize
	\setlength{\tabcolsep}{3pt}
	\begin{tabular}[t]{l|l||rr|rr|rr||rr|rr|rr}
		& & \multicolumn{6}{c||}{Number of column additions} & \multicolumn{6}{c}{Computation time in seconds} \\
		Dataset & Simplices removed & \multicolumn{2}{c|}{\nameshort} & \multicolumn{2}{c|}{PHAT-op (M)} & \multicolumn{2}{c||}{PHAT (M)} & \multicolumn{2}{c|}{\nameshort} & \multicolumn{2}{c|}{PHAT-op} & \multicolumn{2}{c}{PHAT}\\ 
		 & & \texttt{std} & \texttt{twi} & \texttt{std} & \texttt{twi} & \texttt{std} & \texttt{twi} & \texttt{std} & \texttt{twi} & \texttt{std} & \texttt{twi} & \texttt{std} & \texttt{twi} \\
		\hline
        \texttt{vr\_dsphere}      & 183 (0,1,15,59,108) & \textbf{2640}   & \textbf{2508}   & 89.8    & 42.3    & 44.9    & 20.5    & \textbf{13.7}  & 11.9           & 37.5           & 21.8           & 17.6           & \textbf{10.2} \\
        \texttt{vr\_torus}        & 5 (0,0,0,0,5)       & \textbf{28}     & \textbf{28}     & 195.8   & 158.9   & 97.9    & 78.7    & \textbf{6.0}   & \textbf{6.6}   & 70.4           & 66.7           & 36.3           & 33.2 \\
        \texttt{vr\_senate}       & 20 (0,0,0,20)       & \textbf{990}    & \textbf{990}    & 63.6    & 60.4    & 31.8    & 30.1    & 15.3           & 15.4           & 26.9           & 26.0           & \textbf{12.5}  & \textbf{12.4} \\
        \texttt{er\_clique}       & 3 (0,0,0,3)         & \textbf{0}      & \textbf{0}      & 1150.0  & 1141.6  & 576.8   & 570.8   & \textbf{4.7}   & \textbf{5.1}   & 1196.3         & 1132.2         & 429.9          & 394.8 \\
        \texttt{er\_shuffled}     & 25 (0,0,1,24)       & \textbf{7268}   & \textbf{7278}   & 225.1   & 222.7   & 112.6   & 111.3   & \textbf{21.1}  & \textbf{23.6}  & 128.4          & 136.7          & 55.5           & 61.2 \\
        \texttt{alpha\_cube}      & 5 (0,0,0,5)         & \textbf{4}      & \textbf{4}      & 3.3     & 0.3     & 1.7     & 0.1     & \textbf{0.2}   & \textbf{0.1}   & 3.8            & 0.6            & 1.3            & 0.2 \\
        \texttt{alpha\_swissroll} & 9 (0,1,4,4)         & \textbf{14}     & \textbf{18}     & 195.3   & 2.7     & 97.6    & 0.4     & \textbf{10.3}  & 3.0            & 241.3          & 7.1            & 68.2           & \textbf{2.2} \\
        \texttt{cc\_tooth}        & 27 (1,6,12,8)       & \textbf{1228}   & \textbf{138}    & 344.9   & 15.9    & 172.4   & 4.9     & \textbf{70.5}  & 25.5           & 246.7          & 32.9           & 93.5           & \textbf{12.3} \\
        \texttt{bio\_bbmcl6}      & 8 (0,1,6,1)         & \textbf{30}     & \textbf{30}     & 320.9   & 207.9   & 160.5   & 103.8   & \textbf{11.0}  & \textbf{6.3}   & 303.7          & 155.7          & 77.9           & 42.7
     \end{tabular}
	\setlength{\tabcolsep}{6pt}
	\caption{Experiments comparison between \nameshort, PHAT-op, and PHAT, indicating the number of simplices removed (in parentheses in increasing dimension) from each filtration. Number of column additions (in millions ``M'' for PHAT-op and PHAT) and computations time (in seconds, excluding file reading and writing) are given for the standard reduction method and with the twist optimization \cite{twist}, abbreviated \texttt{std} and \texttt{twi}, respectively. The best results are indicated in bold. A visual representation of the testing process and results is given in \cref{fig_testing_environment}.}
\label{table_updates}
\end{table}

\subsubsection{Testing pipeline}
\label{sec_testingpipeline}
The performance of \nameshort~was tested against PHAT-op, the name we give to the PHAT software modified to keep track of the operations matrix $V$ in the decomposition $R=DV$.
This is done so that both programs perform the same task, and the work done is comparable.
The performance of PHAT was also analyzed on the same task, but without requiring to keep track of the column additions, to demonstrate that even if only the barcode is requested (together with the representative cycles of the finite bars), \nameshort~still executes less column additions in a comparable time.
The comparisons among these three programs on the described inputs is given in \cref{table_updates}, and a visual overview of the pipeline and key results is given in \cref{fig_testing_environment}.
Counting the number of column operations was tested separately from timing the computation, as the recording and output the count is itself a time-intensive process.
Read and write times of files were recorded, but are not considered here.
The full code and execution steps are provided in a public repository \cite{our_zenodo}.

\begin{figure}[h]\centering
	\input{figs/testing_environment.tex}
	\caption{
	Steps and outputs of full testing environment.
	The values reported in \cref{table_updates} are from execution of \nameshort and PHAT-op in the ``Testing" stage (right).
	}
	\label{fig_testing_environment}
\end{figure}

\subsubsection{Counting column additions}
The strength of \nameshort~is from the vastly smaller number of column additions that need to be done to update the representative cycles and barcode, as shown in \cref{table_updates}.
Note that in the case of \texttt{er\_clique}, there are no column additions that are done by \nameshort, indicating that the simplices to be removed were not added to other simplices in the reduction, and only their columns need to be deleted and the complex reindexed.
In the same situation PHAT executes more than 500 million column additions.

\subsubsection{Timing algorithm execution}
\label{sec_timing_execution}
As seen in \cref{table_updates}, in almost all cases \nameref{alg_rem} is faster than recomputing from scratch using PHAT, and is always considerably faster than using PHAT-op.
The difference is not only related to the number of simplices removed, as for \texttt{vr\_dsphere} considering the impact of $\approx 180$ simplices is still faster than computing from scratch.
Recall that only the computation time is recorded in \cref{table_updates}, not the file reading in and writing out times.
As demonstrated in \cref{fig_testing_environment} (bottom right), the file sizes may be quite large, so reading and writing them often adds an additional time cost, which skews the results heavily in favor of recomputing from scratch.
However, in the case that the matrices are loaded in memory and are repeatedly being updated, the read and write times will be much less relevant.
The data structure for operations matrices is the same as the one for boundary matrices in PHAT.
These experiments suggest that \nameref{alg_rem} is more efficient in practice under small changes, with room for improvement in memory management and appropriate data structures.

\section{Discussion and future work}
\label{sec_conclusion}

We presented an efficient method to update the barcode and representative cycles of a filtration when simplices are removed.
We proved that our method is asymptotically better than recomputing the barcode from scratch, minimal in the amount of column additions, and we provided constructions on which our method performs asymptotically better than other methods in the literature which can be adapted to produce our output \cite{dey2022updating,vines_and_vineyards}.
Finally, we validated our theoretical results with experiments, demonstrating a time-for-space tradeoff.
We showed that the time necessary to update the barcode and representative cycles with our algorithm is almost always less than recomputing from scratch, especially with a small number of changes, although large input and output files are necessary.

A key improvement necessary for our implementation is decreasing the reading and writing time, as it is now a bottleneck for practical uses, and optimize the look-up of columns in \nameref{alg_rem}, since the number of column additions is already reduced to the theoretical minimum.
In particular, a more efficient method to find all columns with $V[j,i]\neq 0$ for a given $j$ would provide a significant speed up, as currently the algorithm iterates through every column in a given dimension for each simplex to be removed.
Note that \nameref{alg_rem} can also be modified to include several variants of the \nameref{alg_sba} \cite{keepitsparse}.
Finally, we plan to implement a fully dynamic barcode pipeline, merging the existing algorithms for all three type of changes in the filtration (simplex additions, simplex swaps, and simplex removal).
A general algorithm and implementation for all changes exists \cite{luo2023accelerating}, but with more than the minimal number of column additions for simplex removal of \nameref{alg_rem}.

Our experimental comparison with PHAT running in homology is fair for all cases in which the number of elements in the filtration does not grow significantly with the degree, for example with cubical complexes.
However, whenever the data needs to be analyzed using Vietoris--Rips or \v Cech complexes, it is well-known that the \nameref{alg_sba} is much more efficient when run in cohomology with the clearing optimization.
Thus, in this case, we expect our current method to be slower than recomputing from scratch.
Therefore, further work needs to be done to adapt \nameref{alg_rem} to the cohomology case.

Applications of our method, discussed in \cref{sec_motivations_applications}, may require to remove a half or a third of simplices in a given filtration.
In the case that the columns of the simplices to be removed interact with each other and not with the columns of the simplices that remain, one could try to skip even more operations than presented here.
However, implementing this observation is not straightforward and will require care, as some operations affecting only columns-to-be-removed may compound and impact the columns of simplices that stay.
We plan to improve our implementation to take this scenario into account, improving efficiency by considering overall implications of the full set of simplices to be removed.

\section*{Acknowledgements} 
B.G. was partially supported by the Austrian Science Fund (FWF) P 33765-N.
J.L. was partially supported by EPSRC grant EP/P025072/1 ``Topological Analysis of Neural Systems'', and by the Latvian Council of Science (LZP) 1.1.1.9 Research application No 1.1.1.9/LZP/1/24/125 of the Activity ``Post-doctoral Research" ``Efficient topological signatures for representation learning in medical imaging".
J.L. thanks the Austrian Science Fund (FWF) P 33765-N and Prof. Michael Kerber at Graz University of Technology for hosting him during a research visit on this project.
The authors are thankful to Michael Kerber, David Millman, and Felix Ye for useful preliminary discussions, and to anonymous reviewers for constructive feedback.

\renewcommand*{\bibfont}{\small}
\printbibliography
	
\end{document}

%% file: figs/matrix-simplex.tex
\setlength{\arraycolsep}{1pt}
\definecolor{mygray}{gray}{0.4}
\begin{tikzpicture}
\foreach \x\y\n in {0/.2/0, 1.2/1.8/1, 2/.1/2, 2.5/3/7}{
  \coordinate (\n) at (\x,\y);
}
\fill[color0!30] (0)--(1)--(2);
\foreach \x\y\n in {0/1/4, 1/2/5, 1/7/8}{
  \draw (\x) to node [auto] {\n} (\y);
}
\draw (0) to node [auto,swap] {3} (2);
\foreach \x\y\n in {0/.2/0, 1/1.5/1, 2/.1/2, 2.5/3/7}{
  \node[circle,draw,fill=white,minimum size=.5cm] at (\n) {\n};
}

\node at (1.1,.8) {6};
\node (m) at (6.7,1.3) {
$\begin{blockarray}{cccccccccc} 
& {\color{mygray}0} & {\color{mygray}1} & {\color{mygray}2} & {\color{mygray}3} & {\color{mygray}4} & {\color{mygray}5} & {\color{mygray}6} & {\color{mygray}7} & {\color{mygray}8} \\
\begin{block}{c[ccccccccc]}
{\color{mygray}0} & 0 & 0 & 0 & 1 & 1 & 0 & 0 & 0 & 0 \\
{\color{mygray}1} & 0 & 0 & 0 & 0 & 1 & 1 & 0 & 0 & 1 \\ {\color{mygray}2} & 0 & 0 & 0 & 1 & 0 & 1 & 0 & 0 & 0 \\
{\color{mygray}3} & 0 & 0 & 0 & 0 & 0 & 0 & 1 & 0 & 0 \\
{\color{mygray}4} & 0 & 0 & 0 & 0 & 0 & 0 & 1 & 0 & 0 \\
{\color{mygray}5} & 0 & 0 & 0 & 0 & 0 & 0 & 1 & 0 & 0 \\
{\color{mygray}6} & 0 & 0 & 0 & 0 & 0 & 0 & 0 & 0 & 0 \\
{\color{mygray}7} & 0 & 0 & 0 & 0 & 0 & 0 & 0 & 0 & 1 \\
{\color{mygray}8} & 0 & 0 & 0 & 0 & 0 & 0 & 0 & 0 & 0 \\
\end{block}
\end{blockarray}$
};
\begin{scope}[shift={(11,.7)},xscale=1.2]
\draw (-.2,0)--(2.9,0);
\draw[dashed,gray] (-.2,.8)--(2.9,.8);
\foreach \x\n in {0/0,.3/1,.6/2,.9/3,1.2/4,1.5/5,1.8/6,2.1/7,2.4/8,2.7/{\infty}}{
  \draw (\x,-.1)--(\x,.1);
  \node at (\x,-.3) {$\n$};
}
\draw[color0,line width=2pt] (-.02,.2)--(2.72,.2); 
\draw[color0,line width=2pt] (.28,.3)--(1.22,.3); 
\draw[color0,line width=2pt] (.58,.4)--(.92,.4); 
\draw[color0,line width=2pt] (2.18,.5)--(2.42,.5); 
\draw[color0,line width=2pt] (1.48,1.1)--(1.82,1.1); 
\end{scope}
\end{tikzpicture}

%% file: figs/matrix_reduction.tex
\begin{tikzpicture}
\begin{scope}[shift={(-3.5,1.1)}]
	\foreach \r in {1,2,3}{
		\node[circle,draw] (\r) at ($(0,0)+(\r*120-30:.9)$) {$v_\r$};
		\node at ($(0,0)+(270-\r*120:1)$) {$e_\r$};
	}
	\draw (1)--(2);
	\draw[dashed] (1)--(3)--(2);
\end{scope}
\node[anchor=south] at (0,0) {
	\begin{blockarray}{cccc} 
		& $e_1$ & $e_2$ & $e_3$ \\
		\begin{block}{c[ccc]}
			$v_1$ & 1 & 1 & 0 \\
			$v_2$ & 1 & 0 & 1 \\
			$v_3$ & 0 & 1 & 1 \\
		\end{block}
	\end{blockarray}
};
\node (f0) at (.4,0) {(a)};
\node[anchor=south] at (3.5,0) {
	\begin{blockarray}{ccc} 
		$e_1$ & $e_2$ & $e_2+e_3$ \\
		\begin{block}{[ccc]}
			1 & 1 & 1 \\
			1 & 0 & 1 \\
			0 & 1 & 0 \\
		\end{block}
	\end{blockarray}
};
\node (f1) at (3.5,0) {(b)};
\node[anchor=south] at (7.4,0) {
	\begin{blockarray}{ccc} 
		$e_1$ & $e_2$ & $e_1+e_2+e_3$ \\
		\begin{block}{[ccc]}
			1 & 1 & 0 \\
			1 & 0 & 0 \\
			0 & 1 & 0 \\
		\end{block}
	\end{blockarray}
};
\node (f2) at (7.2,0) {(c)};
\end{tikzpicture}

%% file: figs/sirup_worked_example_new_1.tex
\newcommand\cs{.2} 
\newcommand\ts{.5} 
\mathcal{F}\ = \ 
\begin{tikzpicture}[
  baseline=-3pt,
	connector/.style={line width=1pt,rounded corners=8pt},
	nodestyle/.style={circle,draw,fill=white,minimum size=.4cm,scale=.8},
	edgelabel/.style={scale=.8,fill=white,fill opacity=.8,text opacity=1}
]
\node[nodestyle] (8) at (0,0) {8};
\foreach \ang\n\extra in {0/7/.5, 1/0/.5, 2/9/.5, 3/5/0, 4/6/0, 5/2/0, 6/4/0, 7/1/0, 8/3/0}{
	\node[nodestyle] (\n) at ($(0,0)+(\ang*40-40:1.2+\extra)$) {\n};
}
\foreach \n\l in {7/h, 0/a, 9/d, 5/b, 6/e, 2/f, 4/g, 1/j, 3/k}{
	\draw (8) -- (\n);
	\node[edgelabel] at ($(8)!.5!(\n)$) {$\l$};
}
\draw (9) to node[right=3pt,pos=.4,scale=.8] {$c$} (0);
\draw (0) to node[right=3pt,pos=.6,scale=.8] {$i$} (7);
\end{tikzpicture}\ \ ,
\hspace{1cm}
\begin{tikzpicture}[
  baseline=.7cm,
	gridline/.style={white,line width=.2pt}
]
\node at (5.5*\cs,-4*\cs) {$D''$};
\fill[color0!30] (0,0) rectangle (11*\cs,10*\cs);
\fill[color3!30] (2*\cs,0) rectangle (3*\cs,10*\cs);
\foreach \x\y in {1/0, 1/8, 2/5, 2/8, 4/8, 4/9, 5/6, 5/8, 6/2, 6/8, 7/4, 7/8, 8/7, 8/8, 9/0, 9/7, 10/1, 10/8, 11/3, 11/8}{
	\fill[color0] (\x*\cs-\cs,10*\cs-\y*\cs) rectangle (\x*\cs,9*\cs-\y*\cs);
}
\foreach \x\y in {3/0, 3/9}{
	\fill[color3] (\x*\cs-\cs,10*\cs-\y*\cs) rectangle (\x*\cs,9*\cs-\y*\cs);
}
\foreach \x\l in {1/a, 2/b, 3/c, 4/d, 5/e, 6/f, 7/g, 8/h, 9/i, 10/j, 11/k}{ \node[scale=\ts] at (\x*\cs-.5*\cs, 11*\cs) {\vphantom{Ap}$\l$}; }
\foreach \y in {0,...,9}{ \node[scale=\ts] at (-1*\cs, 9.5*\cs-\y*\cs) {$\y$}; }
\foreach \x in {1,...,10}{\draw[gridline] (\x*\cs,0)--(\x*\cs,10*\cs);}
\foreach \y in {1,...,9}{\draw[gridline] (0,\y*\cs)--(11*\cs,\y*\cs);}
\begin{scope}[shift={(14*\cs,-1*\cs)}]
\node at (5.5*\cs,-3*\cs) {$V''$};
\fill[color0!30] (0,0) rectangle (11*\cs,11*\cs);
\fill[color3!30] (2*\cs,0) rectangle (3*\cs,11*\cs);
\fill[color3!30] (0,8*\cs) rectangle (11*\cs,9*\cs);
\foreach \w in {1,2,4,5,6,7,8,9,10,11}{ \fill[color0] (\w*\cs-\cs,12*\cs-\w*\cs) rectangle (\w*\cs,11*\cs-\w*\cs); }
\foreach \a\b in {1/2, 1/4, 1/5, 1/6, 1/7, 1/8, 1/9, 8/9, 1/10, 1/11}{
	\fill[color0] (\b*\cs-\cs,12*\cs-\a*\cs) rectangle (\b*\cs,11*\cs-\a*\cs);
}
\foreach \a\b in {3/3, 3/4}{
	\fill[color3] (\b*\cs-\cs,12*\cs-\a*\cs) rectangle (\b*\cs,11*\cs-\a*\cs);
}
\foreach \w\l in {1/a, 2/b, 3/c, 4/d, 5/e, 6/f, 7/g, 8/h, 9/i, 10/j, 11/k}{
  \node[scale=\ts] at (\w*\cs-.5*\cs, 12*\cs) {\vphantom{Ap}$\l$}; 
  \node[scale=\ts] at (-1*\cs, 11.5*\cs-\w*\cs) {\vphantom{Ap}$\l$}; 
}
\foreach \w in {1,...,10}{
  \draw[gridline] (\w*\cs,0)--(\w*\cs,11*\cs);
  \draw[gridline] (0,\w*\cs)--(11*\cs,\w*\cs);
}
\end{scope}
\begin{scope}[shift={(-16*\cs,0)}]
\node at (5.5*\cs,-4*\cs) {$R''$};
\fill[color0!30] (0,0) rectangle (11*\cs,10*\cs);
\fill[color3!30] (2*\cs,0) rectangle (3*\cs,10*\cs);
\foreach \x\y in {1/0, 1/8, 2/0, 2/5, 5/0, 5/6, 6/0, 6/2, 7/0, 7/4, 8/0, 8/7, 10/0, 10/1, 11/0, 11/3}{
	\fill[color0] (\x*\cs-\cs,10*\cs-\y*\cs) rectangle (\x*\cs,9*\cs-\y*\cs);
}
\foreach \x\y in {3/0, 3/9}{
	\fill[color3] (\x*\cs-\cs,10*\cs-\y*\cs) rectangle (\x*\cs,9*\cs-\y*\cs);
}
\foreach \x\l in {1/a, 2/b, 3/c, 4/d, 5/e, 6/f, 7/g, 8/h, 9/i, 10/j, 11/k}{ \node[scale=\ts] at (\x*\cs-.5*\cs, 11*\cs) {\vphantom{Ap}$\l$}; }
\foreach \y in {0,...,9}{ \node[scale=\ts] at (-1*\cs, 9.5*\cs-\y*\cs) {$\y$}; }
\foreach \x in {1,...,10}{\draw[gridline] (\x*\cs,0)--(\x*\cs,10*\cs);}
\foreach \y in {1,...,9}{\draw[gridline] (0,\y*\cs)--(11*\cs,\y*\cs);}
\end{scope}
\node at (-3*\cs,5*\cs) {$=$};
\end{tikzpicture}\ .

%% file: figs/sirup_worked_example_new_2.tex
\newcommand\cs{.2} 
\newcommand\ts{.5} 
\begin{tikzpicture}[
  baseline=.7cm,
	gridline/.style={white,line width=.2pt}
]
\node at (5*\cs,-3*\cs) {$D'$};
\fill[color0!30] (0,0) rectangle (10*\cs,10*\cs);
\fill[color3!30] (0,0) rectangle (\cs,10*\cs);
\foreach \x\y in {2/5, 2/8, 3/8, 3/9, 4/6, 4/8, 5/2, 5/8, 6/4, 6/8, 7/7, 7/8, 8/0, 8/7, 9/1, 9/8, 10/3, 10/8}{
	\fill[color0] (\x*\cs-\cs,10*\cs-\y*\cs) rectangle (\x*\cs,9*\cs-\y*\cs);
}
\foreach \x\y in {1/0, 1/8}{
	\fill[color3] (\x*\cs-\cs,10*\cs-\y*\cs) rectangle (\x*\cs,9*\cs-\y*\cs);
}
\foreach \x\l in {1/a, 2/b, 3/d, 4/e, 5/f, 6/g, 7/h, 8/i, 9/j, 10/k}{ \node[scale=\ts] at (\x*\cs-.5*\cs, 11*\cs) {\vphantom{Ap}$\l$}; }
\foreach \y in {0,...,9}{ \node[scale=\ts] at (-1*\cs, 9.5*\cs-\y*\cs) {$\y$}; }
\foreach \x in {1,...,9}{\draw[gridline] (\x*\cs,0)--(\x*\cs,10*\cs);}
\foreach \y in {1,...,9}{\draw[gridline] (0,\y*\cs)--(10*\cs,\y*\cs);}
\begin{scope}[shift={(13*\cs,0)}]
\node at (5*\cs,-3*\cs) {$V'$};
\fill[color0!30] (0,0) rectangle (10*\cs,10*\cs);
\fill[color3!30] (0,0) rectangle (\cs,10*\cs);
\fill[color3!30] (0,9*\cs) rectangle (10*\cs,10*\cs);
\foreach \w in {2,...,10}{ \fill[color0] (\w*\cs-\cs,11*\cs-\w*\cs) rectangle (\w*\cs,10*\cs-\w*\cs); }
\foreach \a\b in {7/8}{
	\fill[color0] (\b*\cs-\cs,11*\cs-\a*\cs) rectangle (\b*\cs,10*\cs-\a*\cs);
}
\foreach \a\b in {1/1, 1/2, 1/3, 1/4, 1/5, 1/6, 1/7, 1/8, 1/9, 1/10}{
	\fill[color3] (\b*\cs-\cs,11*\cs-\a*\cs) rectangle (\b*\cs,10*\cs-\a*\cs);
}
\foreach \w\l in {1/a, 2/b, 3/d, 4/e, 5/f, 6/g, 7/h, 8/i, 9/j, 10/k}{
  \node[scale=\ts] at (\w*\cs-.5*\cs, 11*\cs) {\vphantom{Ap}$\l$}; 
  \node[scale=\ts] at (-1*\cs, 10.5*\cs-\w*\cs) {\vphantom{Ap}$\l$}; 
}
\foreach \w in {1,...,9}{
  \draw[gridline] (\w*\cs,0)--(\w*\cs,10*\cs);
  \draw[gridline] (0,\w*\cs)--(10*\cs,\w*\cs);
}
\end{scope}
\begin{scope}[shift={(-15*\cs,0)}]
\node at (5*\cs,-3*\cs) {$R'$};
\fill[color0!30] (0,0) rectangle (10*\cs,10*\cs);
\fill[color3!30] (0,0) rectangle (\cs,10*\cs);
\foreach \x\y in {2/0, 2/5, 3/0, 3/9, 4/0, 4/6, 5/0, 5/2, 6/0, 6/4, 7/0, 7/7, 9/0, 9/1, 10/0, 10/3}{
	\fill[color0] (\x*\cs-\cs,10*\cs-\y*\cs) rectangle (\x*\cs,9*\cs-\y*\cs);
}
\foreach \x\y in {1/0, 1/8}{
	\fill[color3] (\x*\cs-\cs,10*\cs-\y*\cs) rectangle (\x*\cs,9*\cs-\y*\cs);
}
\foreach \x\l in {1/a, 2/b, 3/d, 4/e, 5/f, 6/g, 7/h, 8/i, 9/j, 10/k}{ \node[scale=\ts] at (\x*\cs-.5*\cs, 11*\cs) {\vphantom{Ap}$\l$}; }
\foreach \y in {0,...,9}{ \node[scale=\ts] at (-1*\cs, 9.5*\cs-\y*\cs) {$\y$}; }
\foreach \x in {1,...,9}{\draw[gridline] (\x*\cs,0)--(\x*\cs,10*\cs);}
\foreach \y in {1,...,9}{\draw[gridline] (0,\y*\cs)--(10*\cs,\y*\cs);}
\end{scope}
\node at (-3*\cs,5*\cs) {$=$};
\end{tikzpicture}\ .

%% file: figs/sirup_worked_example_new_3.tex
\newcommand\cs{.2} 
\newcommand\ts{.5} 
\newcommand\css{.05} 
R_N\ =\ 
\begin{tikzpicture}[
  baseline=.7cm,
	gridline/.style={white,line width=.2pt}
]
\fill[color0!30] (0,0) rectangle (9*\cs,10*\cs);
\foreach \x\y in {
  1/5, 2/9, 3/6, 4/2, 5/4, 7/0, 6/7, 8/1, 9/3,
  1/8, 2/8, 3/8, 4/8, 5/8, 6/8, 7/8, 8/8, 9/8
}{
	\fill[color0] (\x*\cs-\cs,10*\cs-\y*\cs) rectangle (\x*\cs,9*\cs-\y*\cs);
}
\foreach \x\l in {1/b, 2/d, 3/e, 4/f, 5/g, 6/h, 7/i, 8/j, 9/k}{ \node[scale=\ts] at (\x*\cs-.5*\cs, 11*\cs) {\vphantom{Ap}$\l$}; }
\foreach \y in {0,...,9}{ \node[scale=\ts] at (-1*\cs, 9.5*\cs-\y*\cs) {$\y$}; }
\foreach \x in {1,...,8}{\draw[gridline] (\x*\cs,0)--(\x*\cs,10*\cs);}
\foreach \y in {1,...,9}{\draw[gridline] (0,\y*\cs)--(9*\cs,\y*\cs);}
\end{tikzpicture}
\ ,\ 
V_N\ =\ 
\begin{tikzpicture}[
  baseline=.7cm,
	gridline/.style={white,line width=.2pt}
]
\fill[color0!30] (0,0) rectangle (9*\cs,9*\cs);
\foreach \w in {1,...,9}{ \fill[color0] (\w*\cs-\cs,10*\cs-\w*\cs) rectangle (\w*\cs,9*\cs-\w*\cs); }
\foreach \a\b in {6/7}{
	\fill[color0] (\b*\cs-\cs,10*\cs-\a*\cs) rectangle (\b*\cs,9*\cs-\a*\cs);
}
\foreach \w\l in {1/b, 2/d, 3/e, 4/f, 5/g, 6/h, 7/i, 8/j, 9/k}{
  \node[scale=\ts] at (\w*\cs-.5*\cs, 10*\cs) {\vphantom{Ap}$\l$}; 
  \node[scale=\ts] at (-1*\cs, 9.5*\cs-\w*\cs) {\vphantom{Ap}$\l$}; 
}
\foreach \w in {1,...,8}{
  \draw[gridline] (\w*\cs,0)--(\w*\cs,9*\cs);
  \draw[gridline] (0,\w*\cs)--(9*\cs,\w*\cs);
}
\end{tikzpicture}
\ ,\ 
R_S\ =\ 
\begin{tikzpicture}[
  baseline=.7cm,
	gridline/.style={white,line width=.2pt}
]
\fill[color0!30] (0,0) rectangle (10*\cs,10*\cs);
\fill[color3!30] (0,0) rectangle (\cs,10*\cs);
\foreach \x\y in {2/0, 2/5, 3/0, 3/9, 4/0, 4/6, 5/0, 5/2, 6/0, 6/4, 7/0, 7/7, 9/0, 9/1, 10/0, 10/3}{
	\fill[color0] (\x*\cs-\cs,10*\cs-\y*\cs) rectangle (\x*\cs,9*\cs-\y*\cs);
}
\foreach \x\y in {1/0, 1/8}{
	\fill[color3] (\x*\cs-\cs,10*\cs-\y*\cs) rectangle (\x*\cs,9*\cs-\y*\cs);
}
\foreach \x\l in {1/a, 2/b, 3/d, 4/e, 5/f, 6/g, 10/k}{ \node[scale=\ts] at (\x*\cs-.5*\cs, 11*\cs) {\vphantom{Ap}$\l$}; }
\foreach \x\l in {7/h, 8/i, 9/j}{ \node[scale=\ts] at (\x*\cs-.5*\cs, 12.5*\cs) {\vphantom{Ap}$\l$}; }
\foreach \y in {0,...,9}{ \node[scale=\ts] at (-1*\cs, 9.5*\cs-\y*\cs) {$\y$}; }
\foreach \x in {1,...,9}{\draw[gridline] (\x*\cs,0)--(\x*\cs,10*\cs);}
\foreach \y in {1,...,9}{\draw[gridline] (0,\y*\cs)--(10*\cs,\y*\cs);}
\foreach \x\y in {1/8, 2/5, 5/2, 8/{-1}}{
	\draw[line width=1pt] (\x*\cs-\cs-\css,10*\cs-\y*\cs+\css) rectangle (\x*\cs+\css,9*\cs-\y*\cs-\css);
}
\end{tikzpicture}
\ .

%% file: figs/sirup_worked_example_new_comparison.tex
\newcommand\cs{.2} 
\newcommand\ts{.5} 
\newcommand\css{.05} 
\newcommand\rowskipper{2.4} 
\newcommand\colskipper{8.5} 
\newcommand\textskipper{15} 
\newcommand\textlinebreak{.1cm} 
\newcommand\naiveR[3]{
	\fill[color0!30] (0,0) rectangle (9*\cs,10*\cs);
  \fill[color2!30] (#1*\cs,0) rectangle (#1*\cs+\cs,10*\cs);
  \foreach \x\y in {#2}{
	  \fill[color0] (\x*\cs-\cs,10*\cs-\y*\cs) rectangle (\x*\cs,9*\cs-\y*\cs);}
  \foreach \x\y in {#3}{
	  \fill[color2] (\x*\cs-\cs,10*\cs-\y*\cs) rectangle (\x*\cs,9*\cs-\y*\cs);}
  \foreach \x\l in {1/b, 2/d, 3/e, 4/f, 5/g, 6/h, 7/i, 8/j, 9/k}{
    \node[scale=\ts] at (\x*\cs-.5*\cs, 10.8*\cs) {\vphantom{Ap}$\l$};}
  \foreach \y in {0,...,9}{\node[scale=\ts] at (-1*\cs, 9.5*\cs-\y*\cs) {$\y$};}
  \foreach \x in {1,...,8}{\draw[gridline] (\x*\cs,0)--(\x*\cs,10*\cs);}
  \foreach \y in {1,...,9}{\draw[gridline] (0,\y*\cs)--(9*\cs,\y*\cs);}}
\newcommand\naiveV[3]{
	\fill[color0!30] (0,0) rectangle (9*\cs,9*\cs);
  \fill[color2!30] (#1*\cs,0) rectangle (#1*\cs+\cs,9*\cs);
  \fill[color2!30] (0,9*\cs-#1*\cs) rectangle (9*\cs,8*\cs-#1*\cs);
  \foreach \w in {1,...,9}{ \fill[color0] (\w*\cs-\cs,10*\cs-\w*\cs) rectangle (\w*\cs,9*\cs-\w*\cs); }
  \foreach \a\b in {#2}{\fill[color0] (\b*\cs-\cs,10*\cs-\a*\cs) rectangle (\b*\cs,9*\cs-\a*\cs);}
  \foreach \a\b in {#3}{\fill[color2] (\b*\cs-\cs,10*\cs-\a*\cs) rectangle (\b*\cs,9*\cs-\a*\cs);}
  \foreach \w\l in {1/b, 2/d, 3/e, 4/f, 5/g, 6/h, 7/i, 8/j, 9/k}{
    \node[scale=\ts] at (\w*\cs-.5*\cs, 9.8*\cs) {\vphantom{Ap}$\l$}; 
    \node[scale=\ts] at (-1*\cs, 9.5*\cs-\w*\cs) {\vphantom{Ap}$\l$};}
  \foreach \w in {1,...,8}{
    \draw[gridline] (\w*\cs,0)--(\w*\cs,9*\cs);
    \draw[gridline] (0,\w*\cs)--(9*\cs,\w*\cs);}}
\newcommand\sirupR[6]{
	\fill[color0!30] (0,0) rectangle (10*\cs,10*\cs);
  \fill[color2!30] (#1*\cs,0) rectangle (#1*\cs+\cs,10*\cs);
  \fill[color3!30] (0,0) rectangle (\cs,10*\cs);
  \foreach \x\y in {#2}{\fill[color0] (\x*\cs-\cs,10*\cs-\y*\cs) rectangle (\x*\cs,9*\cs-\y*\cs);}
  \foreach \x\y in {#3}{\fill[color2] (\x*\cs-\cs,10*\cs-\y*\cs) rectangle (\x*\cs,9*\cs-\y*\cs);}
  \foreach \x\y in {1/0, 1/8}{\fill[color3] (\x*\cs-\cs,10*\cs-\y*\cs) rectangle (\x*\cs,9*\cs-\y*\cs);}
  \foreach \x\l in {1/a, 2/b, 3/d, 4/e, 5/f, 6/g, 7/h, 9/j, 10/k}{
    \node[scale=\ts] at (\x*\cs-.5*\cs, 10.8*\cs) {\vphantom{Ap}$\l$}; }
  \foreach \x\l in {8/i}{
    \node[scale=\ts] at (\x*\cs-.5*\cs, 10.7*\cs) {\vphantom{Ap}$\l$}; }
  \fill[white] (#5*\cs-\cs,10.1*\cs) rectangle (#6*\cs,11.5*\cs);
  \foreach \y in {0,...,9}{ \node[scale=\ts] at (-1*\cs, 9.5*\cs-\y*\cs) {$\y$}; }
  \foreach \x in {1,...,9}{\draw[gridline] (\x*\cs,0)--(\x*\cs,10*\cs);}
  \foreach \y in {1,...,9}{\draw[gridline] (0,\y*\cs)--(10*\cs,\y*\cs);}
  \foreach \x\y in {8/{-1}, 5/2, 2/5, 1/8}{
    \draw[line width=1pt,opacity=.2] (\x*\cs-\cs-\css,10*\cs-\y*\cs+\css) rectangle (\x*\cs+\css,9*\cs-\y*\cs-\css);}
  \foreach \x\y in {#4}{}}
\newcommand\sirupV[4]{
	\fill[color0!30] (0,0) rectangle (10*\cs,10*\cs);
  \fill[color3!30] (0,0) rectangle (\cs,10*\cs);
  \fill[color3!30] (0,9*\cs) rectangle (10*\cs,10*\cs);
  \fill[color2!30] (#1*\cs,0) rectangle (#1*\cs+\cs,10*\cs);
  \fill[color2!30] (0,9*\cs-#1*\cs) rectangle (10*\cs,10*\cs-#1*\cs);
  \foreach \w in {2,...,10}{ \fill[color0] (\w*\cs-\cs,11*\cs-\w*\cs) rectangle (\w*\cs,10*\cs-\w*\cs); }
  \foreach \a\b in {#2}{\fill[color0] (\b*\cs-\cs,11*\cs-\a*\cs) rectangle (\b*\cs,10*\cs-\a*\cs);}
  \foreach \a\b in {#3}{\fill[color2] (\b*\cs-\cs,11*\cs-\a*\cs) rectangle (\b*\cs,10*\cs-\a*\cs);}
  \foreach \a\b in {#4}{\fill[color3] (\b*\cs-\cs,11*\cs-\a*\cs) rectangle (\b*\cs,10*\cs-\a*\cs);}
  \foreach \w\l in {1/a, 2/b, 3/d, 4/e, 5/f, 6/g, 7/h, 8/i, 9/j, 10/k}{
    \node[scale=\ts] at (\w*\cs-.5*\cs, 10.8*\cs) {\vphantom{Ap}$\l$}; 
    \node[scale=\ts] at (-1*\cs, 10.5*\cs-\w*\cs) {\vphantom{Ap}$\l$};}
  \foreach \w in {1,...,9}{
    \draw[gridline] (\w*\cs,0)--(\w*\cs,10*\cs);
    \draw[gridline] (0,\w*\cs)--(10*\cs,\w*\cs);}}
\newcommand\connectorbelow[4]{
  \draw[{Circle[width=2pt,length=2pt]}-{Circle[width=2pt,length=2pt]},line width=.5pt,shorten <=-1pt,shorten >=-1pt]
  (#1*\cs-.5*\cs,10*\cs-#2*\cs-.5*\cs) |- (#3*\cs-.5*\cs,10*\cs-#4*\cs-.5*\cs);}
\newcommand\connectorabove[4]{
  \draw[{Circle[width=2pt,length=2pt]}-{Circle[width=2pt,length=2pt]},line width=.5pt,shorten <=-1pt,shorten >=-1pt]
  (#1*\cs-.5*\cs,10*\cs-#2*\cs-.5*\cs) -| (#3*\cs-.5*\cs,10*\cs-#4*\cs-.5*\cs);}
\newcommand\simpleconnector[3]{
  \draw[{Circle[width=2pt,length=2pt]}-{Straight Barb[width=3pt,length=1.5pt]},line width=.5pt, shorten <=-1pt, shorten >=2pt]
  (#1*\cs-.5*\cs,10*\cs-#3*\cs-.5*\cs) -- (#2*\cs,10*\cs-#3*\cs-.5*\cs);}
\begin{tikzpicture}[gridline/.style={white,line width=.2pt}]
\begin{scope}[shift={(0,8*\rowskipper)}]
\node[anchor=west] (text) at (-\textskipper*\cs,10*\cs) {\vphantom{Ap}1. Column $b$:};
\node[anchor=north west,xshift=.3cm,yshift=\textlinebreak] at (text.south west) {\vphantom{Ap}-};
\naiveR{0}
  {2/8, 2/9, 3/6, 3/8, 4/2, 4/8, 5/4, 5/8, 6/7, 6/8, 7/0, 7/8, 8/1, 8/8, 9/3, 9/8}
  {1/5, 1/8}
\begin{scope}[shift={(12*\cs,\cs)}]
\naiveV{0}{6/7}{1/1}
\end{scope}\end{scope}
\begin{scope}[shift={(0,7*\rowskipper)}]
\node[anchor=west] (text) at (-\textskipper*\cs,10*\cs) {\vphantom{Ap}2. Column $d$:};
\node[anchor=north west,xshift=.3cm,yshift=\textlinebreak] at (text.south west) {\vphantom{Ap}-};
\naiveR{1}
  {1/5, 1/8, 3/6, 3/8, 4/2, 4/8, 5/4, 5/8, 6/7, 6/8, 7/0, 7/8, 8/1, 8/8, 9/3, 9/8}
  {2/8, 2/9}
\begin{scope}[shift={(12*\cs,\cs)}]
\naiveV{1}{6/7}{2/2}
\end{scope}\end{scope}
\begin{scope}[shift={(0,6*\rowskipper)}]
\node[anchor=west] (text) at (-\textskipper*\cs,10*\cs) {\vphantom{Ap}3. Column $e$:};
\node[anchor=north west,xshift=.3cm,yshift=\textlinebreak] at (text.south west) {\vphantom{Ap}add $b$};
\naiveR{2}
  {1/5, 1/8, 2/8, 2/9, 4/2, 4/8, 5/4, 5/8, 6/7, 6/8, 7/0, 7/8, 8/1, 8/8, 9/3, 9/8}
  {3/5, 3/6}
\simpleconnector{1}{3}{8}
\begin{scope}[shift={(12*\cs,\cs)}]
\naiveV{2}{6/7}{1/3, 3/3}
\end{scope}\end{scope}
\begin{scope}[shift={(0,5*\rowskipper)}]
\node[anchor=west] (text) at (-\textskipper*\cs,10*\cs) {\vphantom{Ap}4. Column $f$:};
\node[anchor=north west,xshift=.3cm,yshift=\textlinebreak] at (text.south west) {\vphantom{Ap}add $b$};
\naiveR{3}
  {1/5, 1/8, 2/8, 2/9, 3/5, 3/6, 5/4, 5/8, 6/7, 6/8, 7/0, 7/8, 8/1, 8/8, 9/3, 9/8}
  {4/2, 4/5}
\simpleconnector{1}{4}{8}
\begin{scope}[shift={(12*\cs,\cs)}]
\naiveV{3}{1/3, 6/7}{1/4, 4/4}
\end{scope}\end{scope}
\begin{scope}[shift={(0,4*\rowskipper)}]
\node[anchor=west] (text1) at (-\textskipper*\cs,10*\cs) {\vphantom{Ap}5. Column $g$:};
\node[anchor=north west,xshift=.3cm,yshift=\textlinebreak] (text2) at (text1.south west) {\vphantom{Ap}add $b$};
\node[anchor=north west,yshift=\textlinebreak] (text3) at (text2.south west) {\vphantom{Ap}add $f$};
\naiveR{4}
  {1/5, 1/8, 2/8, 2/9, 3/5, 3/6, 4/2, 4/5, 6/7, 6/8, 7/0, 7/8, 8/1, 8/8, 9/3, 9/8}
  {5/2, 5/4}
\simpleconnector{1}{5}{8}
\simpleconnector{4}{5}{5}
\begin{scope}[shift={(12*\cs,\cs)}]
\naiveV{4}{1/3, 1/4, 6/7}{4/5, 5/5}
\end{scope}\end{scope}
\begin{scope}[shift={(0,3*\rowskipper)}]
\node[anchor=west] (text1) at (-\textskipper*\cs,10*\cs) {\vphantom{Ap}6. Column $h$:};
\node[anchor=north west,xshift=.3cm,yshift=\textlinebreak] (text2) at (text1.south west) {\vphantom{Ap}add $b$};
\naiveR{5}
  {1/5, 1/8, 2/8, 2/9, 3/5, 3/6, 4/2, 4/5, 5/2, 5/4, 7/0, 7/8, 8/1, 8/8, 9/3, 9/8}
  {6/5, 6/7}
\simpleconnector{1}{6}{8}
\begin{scope}[shift={(12*\cs,\cs)}]
\naiveV{5}{1/3, 1/4, 4/5}{1/6, 6/6, 6/7}
\end{scope}\end{scope}
\begin{scope}[shift={(0,2*\rowskipper)}]
\node[anchor=west] (text1) at (-\textskipper*\cs,10*\cs) {\vphantom{Ap}7. Column $i$:};
\node[anchor=north west,xshift=.3cm,yshift=\textlinebreak] (text2) at (text1.south west) {\vphantom{Ap}add $b$};
\node[anchor=north west,yshift=\textlinebreak] (text3) at (text2.south west) {\vphantom{Ap}add $f$};
\naiveR{6}
  {1/5, 1/8, 2/8, 2/9, 3/5, 3/6, 4/2, 4/5, 5/2, 5/4, 6/5, 6/7, 8/1, 8/8, 9/3, 9/8}
  {7/0, 7/2}
\simpleconnector{1}{7}{8}
\simpleconnector{4}{7}{5}
\begin{scope}[shift={(12*\cs,\cs)}]
\naiveV{6}{1/3, 1/4, 4/5, 1/6}{4/7, 6/7, 7/7}
\end{scope}\end{scope}
\begin{scope}[shift={(0,1*\rowskipper)}]
\node[anchor=west] (text1) at (-\textskipper*\cs,10*\cs) {\vphantom{Ap}8. Column $j$:};
\node[anchor=north west,xshift=.3cm,yshift=\textlinebreak] (text2) at (text1.south west) {\vphantom{Ap}add $b$};
\node[anchor=north west,yshift=\textlinebreak] (text3) at (text2.south west) {\vphantom{Ap}add $f$};
\node[anchor=north west,yshift=\textlinebreak] (text4) at (text3.south west) {\vphantom{Ap}add $i$};
\naiveR{7}
  {1/5, 1/8, 2/8, 2/9, 3/5, 3/6, 4/2, 4/5, 5/2, 5/4, 6/5, 6/7, 7/0, 7/2, 9/3, 9/8}
  {8/0, 8/1}
\simpleconnector{1}{8}{8}
\simpleconnector{4}{8}{5}
\simpleconnector{7}{8}{2}
\begin{scope}[shift={(12*\cs,\cs)}]
\naiveV{7}{1/3, 1/4, 4/5, 1/6, 4/7, 6/7}{6/8, 7/8, 8/8}
\end{scope}\end{scope}
\begin{scope}[shift={(0,0*\rowskipper)}]
\node[anchor=west] (text1) at (-\textskipper*\cs,10*\cs) {\vphantom{Ap}9. Column $k$:};
\node[anchor=north west,xshift=.3cm,yshift=\textlinebreak] (text2) at (text1.south west) {\vphantom{Ap}add $b$};
\node[anchor=north west,yshift=\textlinebreak] (text3) at (text2.south west) {\vphantom{Ap}add $f$};
\naiveR{8}
  {1/5, 1/8, 2/8, 2/9, 3/5, 3/6, 4/2, 4/5, 5/2, 5/4, 6/5, 6/7, 7/0, 7/2, 8/0, 8/1}
  {9/2, 9/3}
\simpleconnector{1}{9}{8}
\simpleconnector{4}{9}{5}
\begin{scope}[shift={(12*\cs,\cs)}]
\naiveV{8}{1/3, 1/4, 4/5, 1/6, 4/7, 6/7, 6/8, 7/8}{4/9, 9/9}
\end{scope}\end{scope}
\begin{scope}[shift={(\colskipper, 8*\rowskipper)}]
\node[anchor=west] (text1) at (-1*\textskipper*\cs,10*\cs) {\vphantom{Ap}1. Column $k$:};
\node[anchor=north west,xshift=.3cm,yshift=\textlinebreak,gray] (text2) at (text1.south west) {\vphantom{Ap}(not in $B$)};
\node[anchor=north west,yshift=\textlinebreak] (text3) at (text2.south west) {\vphantom{Ap}add $f$};
\sirupR{9}
  {2/0, 2/5, 3/0, 3/9, 4/0, 4/6, 5/0, 5/2, 6/0, 6/4, 7/0, 7/7, 9/0, 9/1}
	{10/2, 10/3}{5/2}{0}{0}
\simpleconnector{5}{10}{2}
\begin{scope}[shift={(13*\cs,0)}]
\sirupV{9}{7/8}{5/10, 10/10}{1/1, 1/2, 1/3, 1/4, 1/5, 1/6, 1/7, 1/8, 1/9}
\end{scope}\end{scope}
\begin{scope}[shift={(\colskipper, 7*\rowskipper)}]
\node[anchor=west] (text1) at (-1*\textskipper*\cs,10*\cs) {\vphantom{Ap}2. Column $j$:};
\node[anchor=north west,xshift=.3cm,yshift=\textlinebreak,gray] (text2) at (text1.south west) {\vphantom{Ap}(not in $B$)};
\node[anchor=north west,yshift=\textlinebreak] (text3) at (text2.south west) {\vphantom{Ap}add $i$};
\sirupR{8}
  {2/0, 2/5, 3/0, 3/9, 4/0, 4/6, 5/0, 5/2, 6/0, 6/4, 7/0, 7/7, 10/2, 10/3}
	{9/0, 9/1}{8/{-1}}{8}{9}
\simpleconnector{8}{9}{-1}
\begin{scope}[shift={(13*\cs,0)}]
\sirupV{8}{7/8, 5/10, 10/10}{7/9, 8/9,9/9}{1/1, 1/2, 1/3, 1/4, 1/5, 1/6, 1/7, 1/8}
\end{scope}\end{scope}
\begin{scope}[shift={(\colskipper, 6*\rowskipper)}]
\node[anchor=west] (text1) at (-1*\textskipper*\cs,10*\cs) {\vphantom{Ap}3. Column $i$:};
\node[anchor=north west,xshift=.3cm,yshift=\textlinebreak,gray] (text2) at (text1.south west) {\vphantom{Ap}(in $B$)};
\node[anchor=north west,yshift=\textlinebreak] (text3) at (text2.south west) {\vphantom{Ap}add $f$};
\sirupR{7}
  {2/0, 2/5, 3/0, 3/9, 4/0, 4/6, 5/0, 5/2, 6/0, 6/4, 7/0, 7/7, 9/0, 9/1, 10/2, 10/3}
	{8/0, 8/2}{5/2, 8/{-1}}{0}{0}
\simpleconnector{5}{8}{2}
\begin{scope}[shift={(13*\cs,0)}]
\sirupV{7}{7/9, 5/10}{5/8, 7/8, 8/8, 8/9}{1/1, 1/2, 1/3, 1/4, 1/5, 1/6, 1/7}
\end{scope}\end{scope}
\begin{scope}[shift={(\colskipper, 5*\rowskipper)}]
\node[anchor=west] (text1) at (-1*\textskipper*\cs,10*\cs) {\vphantom{Ap}4. Column $h$:};
\node[anchor=north west,xshift=.3cm,yshift=\textlinebreak,gray] (text2) at (text1.south west) {\vphantom{Ap}(not in $B$)};
\node[anchor=north west,yshift=\textlinebreak] (text3) at (text2.south west) {\vphantom{Ap}add $b$};
\sirupR{6}
  {2/0, 2/5, 3/0, 3/9, 4/0, 4/6, 5/0, 5/2, 6/0, 6/4, 8/0, 8/2, 9/0, 9/1, 10/2, 10/3}
	{7/5, 7/7}{2/5}{0}{0}
\simpleconnector{2}{7}{5}
\begin{scope}[shift={(13*\cs,0)}]
\sirupV{6}{5/10, 5/8, 8/9}{2/7, 7/8, 7/9, 7/7}{1/1, 1/2, 1/3, 1/4, 1/5, 1/6}
\end{scope}\end{scope}
\begin{scope}[shift={(\colskipper, 4*\rowskipper)}]
\node[anchor=west] (text1) at (-1*\textskipper*\cs,10*\cs) {\vphantom{Ap}5. Column $g$:};
\node[anchor=north west,xshift=.3cm,yshift=\textlinebreak,gray] (text2) at (text1.south west) {\vphantom{Ap}(not in $B$)};
\node[anchor=north west,yshift=\textlinebreak] (text3) at (text2.south west) {\vphantom{Ap}add $f$};
\sirupR{5}
  {2/0, 2/5, 3/0, 3/9, 4/0, 4/6, 5/0, 5/2, 7/5, 7/7, 8/0, 8/2, 9/0, 9/1, 10/2, 10/3}
	{6/2, 6/4}{5/2}{0}{0}
\simpleconnector{5}{6}{2}
\begin{scope}[shift={(13*\cs,0)}]
\sirupV{5}{5/10, 5/8, 8/9, 2/7, 7/8, 7/9}{5/6, 6/6}{1/1, 1/2, 1/3, 1/4, 1/5}
\end{scope}\end{scope}
\begin{scope}[shift={(\colskipper, 3*\rowskipper)}]
\node[anchor=west] (text1) at (-1*\textskipper*\cs,10*\cs) {\vphantom{Ap}6. Column $f$:};
\node[anchor=north west,xshift=.3cm,yshift=\textlinebreak,gray] (text2) at (text1.south west) {\vphantom{Ap}(in $B$)};
\node[anchor=north west,yshift=\textlinebreak] (text3) at (text2.south west) {\vphantom{Ap}add $b$};
\sirupR{4}
  {2/0, 2/5, 3/0, 3/9, 4/0, 4/6, 6/2, 6/4, 7/5, 7/7, 8/0, 8/2, 9/0, 9/1, 10/2, 10/3}
	{5/2, 5/5}{2/5, 5/2}{0}{0}
\simpleconnector{2}{5}{5}
\begin{scope}[shift={(13*\cs,0)}]
\sirupV{4}{8/9, 2/7, 7/8, 7/9}{2/5, 5/10, 5/8, 5/6, 5/5}{1/1, 1/2, 1/3, 1/4}
\end{scope}\end{scope}
\begin{scope}[shift={(\colskipper, 2*\rowskipper)}]
\node[anchor=west] (text1) at (-1*\textskipper*\cs,10*\cs) {\vphantom{Ap}7. Column $e$:};
\node[anchor=north west,xshift=.3cm,yshift=\textlinebreak,gray] (text2) at (text1.south west) {\vphantom{Ap}(not in $B$)};
\node[anchor=north west,yshift=\textlinebreak] (text3) at (text2.south west) {\vphantom{Ap}add $b$};
\sirupR{3}
  {2/0, 2/5, 3/0, 3/9, 5/2, 5/5, 6/2, 6/4, 7/5, 7/7, 8/0, 8/2, 9/0, 9/1, 10/2, 10/3}
	{4/5, 4/6}{2/5}{0}{0}
\simpleconnector{2}{4}{5}
\begin{scope}[shift={(13*\cs,0)}]
\sirupV{3}{8/9, 2/7, 7/8, 7/9, 5/6, 2/5, 5/10, 5/8}{2/4, 4/4}{1/1, 1/2, 1/3}
\end{scope}\end{scope}
\begin{scope}[shift={(\colskipper, 1*\rowskipper)}]
\node[anchor=west] (text1) at (-1*\textskipper*\cs,10*\cs) {\vphantom{Ap}8. Column $d$:};
\node[anchor=north west,xshift=.3cm,yshift=\textlinebreak,gray] (text2) at (text1.south west) {\vphantom{Ap}(not in $B$)};
\node[anchor=north west,yshift=\textlinebreak] (text3) at (text2.south west) {\vphantom{Ap}add $a$};
\sirupR{2}
  {2/0, 2/5, 4/5, 4/6, 5/2, 5/5, 6/2, 6/4, 7/5, 7/7, 8/0, 8/2, 9/0, 9/1, 10/2, 10/3}
	{3/8, 3/9}{1/8}{0}{0}
\simpleconnector{1}{3}{8}
\begin{scope}[shift={(13*\cs,0)}]
\sirupV{2}{8/9, 2/7, 7/8, 7/9, 5/6, 2/5, 5/10, 5/8, 2/4}{3/3}{1/1, 1/2}
\end{scope}\end{scope}
\begin{scope}[shift={(\colskipper, 0*\rowskipper)}]
\node[anchor=west] (text1) at (-1*\textskipper*\cs,10*\cs) {\vphantom{Ap}9. Column $b$:};
\node[anchor=north west,xshift=.3cm,yshift=\textlinebreak,gray] (text2) at (text1.south west) {\vphantom{Ap}(in $B$)};
\node[anchor=north west,yshift=\textlinebreak] (text3) at (text2.south west) {\vphantom{Ap}add $a$};
\sirupR{1}
  {3/8, 3/9, 4/5, 4/6, 5/2, 5/5, 6/2, 6/4, 7/5, 7/7, 8/0, 8/2, 9/0, 9/1, 10/2, 10/3}
	{2/5, 2/8}{1/8, 2/5}{0}{0}
\simpleconnector{1}{2}{8}
\begin{scope}[shift={(13*\cs,0)}]
\sirupV{1}{8/9, 7/8, 7/9, 5/6, 5/10, 5/8}{2/7, 2/5, 2/4, 2/2}{1/1}
\end{scope}\end{scope}
\end{tikzpicture}

%% file: figs/lemma_bpivots.tex
\newcommand\cs{.2} 
\begin{tikzpicture}[
	baseline=2*\cs cm,
	gridline/.style={white,line width=.2pt}
]
\fill[color0!30] (0,0) rectangle (5*\cs,5*\cs);
\fill[color3!30] (0,0) rectangle (\cs,5*\cs);
\fill[color3] (0,0) rectangle (\cs,\cs);
\foreach \x\y in {1/1, 2/2, 3/3, 4/4}{
	\fill[color0] (\x*\cs,\y*\cs) rectangle (\x*\cs+\cs,\y*\cs+\cs);
}
\foreach \xy in {1,...,4}{\
	\draw[gridline] (\xy*\cs,0)--(\xy*\cs,5*\cs);
	\draw[gridline] (0,\xy*\cs)--(5*\cs,\xy*\cs);
}
\node at (2.5*\cs,-2*\cs) {$R'$};
\begin{scope}[shift={(6*\cs,0)}]
\fill[color0!30] (0,0) rectangle (5*\cs,5*\cs);
\fill[color3!30] (0,0) rectangle (\cs,5*\cs);
\fill[color3!30] (0,5*\cs) rectangle (5*\cs,4*\cs);
\fill[color3] (0,4*\cs) rectangle (\cs,5*\cs);
\foreach \x\y in {1/3, 2/2, 3/1, 4/0}{
	\fill[color0] (\x*\cs,\y*\cs) rectangle (\x*\cs+\cs,\y*\cs+\cs);
}
\foreach \x\y in {1/4, 2/4, 3/4, 4/4}{
	\fill[color3] (\x*\cs,\y*\cs) rectangle (\x*\cs+\cs,\y*\cs+\cs);
}
\foreach \xy in {1,...,4}{\
	\draw[gridline] (\xy*\cs,0)--(\xy*\cs,5*\cs);
	\draw[gridline] (0,\xy*\cs)--(5*\cs,\xy*\cs);
}
\node at (2.5*\cs,-2*\cs) {$V'$};
\end{scope}
\end{tikzpicture}
\ \ \ \text{becomes}\ \ \ 
\begin{tikzpicture}[
	baseline=2*\cs cm,
	gridline/.style={white,line width=.2pt}
]
\fill[color0!30] (0,0) rectangle (5*\cs,5*\cs);
\fill[color3!30] (0,0) rectangle (\cs,5*\cs);
\fill[color3] (0,0) rectangle (\cs,\cs);
\foreach \x\y in {1/1, 1/0, 2/2, 2/1, 3/3, 3/2, 4/4, 4/3}{
	\fill[color0] (\x*\cs,\y*\cs) rectangle (\x*\cs+\cs,\y*\cs+\cs);
}
\foreach \xy in {1,...,4}{\
	\draw[gridline] (\xy*\cs,0)--(\xy*\cs,5*\cs);
	\draw[gridline] (0,\xy*\cs)--(5*\cs,\xy*\cs);
}
\node at (2.5*\cs,-2*\cs) {$R$};
\begin{scope}[shift={(6*\cs,0)}]
\fill[color0!30] (0,0) rectangle (5*\cs,5*\cs);
\fill[color3!30] (0,0) rectangle (\cs,5*\cs);
\fill[color3!30] (0,5*\cs) rectangle (5*\cs,4*\cs);
\fill[color3] (0,4*\cs) rectangle (\cs,5*\cs);
\foreach \x\y in {1/3, 2/3, 2/2, 3/2, 3/1, 4/1, 4/0}{
	\fill[color0] (\x*\cs,\y*\cs) rectangle (\x*\cs+\cs,\y*\cs+\cs);
}
\foreach \xy in {1,...,4}{\
	\draw[gridline] (\xy*\cs,0)--(\xy*\cs,5*\cs);
	\draw[gridline] (0,\xy*\cs)--(5*\cs,\xy*\cs);
}
\node at (2.5*\cs,-2*\cs) {$V$};
\end{scope}
\end{tikzpicture}
\ .

%% file: figs/lemma_apivots.tex
\newcommand\cs{.2} 
\begin{tikzpicture}[
	baseline=2*\cs cm,
	gridline/.style={white,line width=.2pt}
]
\fill[color0!30] (0,0) rectangle (5*\cs,5*\cs);
\fill[color3!30] (0,0) rectangle (\cs,5*\cs);
\fill[color3] (0,4*\cs) rectangle (\cs,5*\cs);
\foreach \x\y in {1/3, 2/2, 3/1, 4/0}{
	\fill[color0] (\x*\cs,\y*\cs) rectangle (\x*\cs+\cs,\y*\cs+\cs);
}
\foreach \xy in {1,...,4}{\
	\draw[gridline] (\xy*\cs,0)--(\xy*\cs,5*\cs);
	\draw[gridline] (0,\xy*\cs)--(5*\cs,\xy*\cs);
}
\node at (2.5*\cs,-2*\cs) {$R'$};
\begin{scope}[shift={(6*\cs,0)}]
\fill[color0!30] (0,0) rectangle (5*\cs,5*\cs);
\fill[color3!30] (0,0) rectangle (\cs,5*\cs);
\fill[color3!30] (0,5*\cs) rectangle (5*\cs,4*\cs);
\fill[color3] (0,4*\cs) rectangle (\cs,5*\cs);
\foreach \x\y in {1/3, 2/2, 3/1, 4/0}{
	\fill[color0] (\x*\cs,\y*\cs) rectangle (\x*\cs+\cs,\y*\cs+\cs);
}
\foreach \x\y in {1/4, 2/4, 3/4, 4/4}{
	\fill[color3] (\x*\cs,\y*\cs) rectangle (\x*\cs+\cs,\y*\cs+\cs);
}
\foreach \xy in {1,...,4}{\
	\draw[gridline] (\xy*\cs,0)--(\xy*\cs,5*\cs);
	\draw[gridline] (0,\xy*\cs)--(5*\cs,\xy*\cs);
}
\node at (2.5*\cs,-2*\cs) {$V'$};
\end{scope}
\end{tikzpicture}
\ \ \ \text{becomes}\ \ \ 
\begin{tikzpicture}[
	baseline=2*\cs cm,
	gridline/.style={white,line width=.2pt}
]
\fill[color0!30] (0,0) rectangle (5*\cs,5*\cs);
\fill[color3!30] (0,0) rectangle (\cs,5*\cs);
\fill[color3] (0,4*\cs) rectangle (\cs,5*\cs);
\foreach \x\y in {1/4, 1/3, 2/4, 2/2, 3/4, 3/1, 4/4, 4/0}{
	\fill[color0] (\x*\cs,\y*\cs) rectangle (\x*\cs+\cs,\y*\cs+\cs);
}
\foreach \xy in {1,...,4}{\
	\draw[gridline] (\xy*\cs,0)--(\xy*\cs,5*\cs);
	\draw[gridline] (0,\xy*\cs)--(5*\cs,\xy*\cs);
}
\node at (2.5*\cs,-2*\cs) {$R$};
\begin{scope}[shift={(6*\cs,0)}]
\fill[color0!30] (0,0) rectangle (5*\cs,5*\cs);
\fill[color3!30] (0,0) rectangle (\cs,5*\cs);
\fill[color3!30] (0,5*\cs) rectangle (5*\cs,4*\cs);
\fill[color3] (0,4*\cs) rectangle (\cs,5*\cs);
\foreach \x\y in {1/3, 2/2, 3/1, 4/0}{
	\fill[color0] (\x*\cs,\y*\cs) rectangle (\x*\cs+\cs,\y*\cs+\cs);
}
\foreach \xy in {1,...,4}{\
	\draw[gridline] (\xy*\cs,0)--(\xy*\cs,5*\cs);
	\draw[gridline] (0,\xy*\cs)--(5*\cs,\xy*\cs);
}
\node at (2.5*\cs,-2*\cs) {$V$};
\end{scope}
\end{tikzpicture}
\ .

%% file: figs/clearing.tex
\newcommand\cs{.24} 
\newcommand\ylabel{-.7} 
\begin{tikzpicture}[
	gridline/.style={white,line width=.2pt},
	arrdir/.style={-{Straight Barb[width=6pt,length=3pt]}}
]
\foreach \x\y\n in {0/.2/0, 1.2/1.8/1, 2/.1/2}{
  \coordinate (\n) at (.9*\x,.9*\y);
}
\fill[color0!30] (0)--(1)--(2);
\foreach \x\y\n in {0/1/4, 1/2/5}{
  \draw (\x) to node [auto] {\n} (\y);
}
\draw (0) to node [auto,swap] {3} (2);
\foreach \x\y\n in {0/.2/0, 1/1.5/1, 2/.1/2}{
  \node[circle,draw,fill=white,minimum size=.5cm] at (\n) {\n};
}
\node at (1,.7) {6};
\node at (1,\ylabel) {$K$};
\begin{scope}[shift={(2.5,0)}]
\fill[color0!30] (0,0) rectangle (7*\cs,7*\cs);
\foreach \x\y in {3/0, 3/2, 4/0, 4/1, 5/1, 5/2, 6/3, 6/4, 6/5}{
  \fill[color0] (\x*\cs,7*\cs-\y*\cs) rectangle (\x*\cs+\cs,6*\cs-\y*\cs);
}
\foreach \xy in {1,...,6}{
  \draw[gridline] (\xy*\cs,0)--(\xy*\cs,7*\cs);
  \draw[gridline] (0,\xy*\cs)--(7*\cs,\xy*\cs);
 }
\node at (3.5*\cs,\ylabel) {$D$};
\end{scope}
\begin{scope}[shift={(5.8,0)}]
\fill[color0!30] (0,0) rectangle (7*\cs,7*\cs);
\fill[color2!30] (3*\cs,0) rectangle (5*\cs,7*\cs);
\foreach \x\y in {6/3, 6/4, 6/5}{
  \fill[color0] (\x*\cs,7*\cs-\y*\cs) rectangle (\x*\cs+\cs,6*\cs-\y*\cs);
}
\foreach \x\y in {3/0, 3/2, 4/0, 4/1}{
  \fill[color2] (\x*\cs,7*\cs-\y*\cs) rectangle (\x*\cs+\cs,6*\cs-\y*\cs);
}
\foreach \xy in {1,...,6}{
  \draw[gridline] (\xy*\cs,0)--(\xy*\cs,7*\cs);
  \draw[gridline] (0,\xy*\cs)--(7*\cs,\xy*\cs);
 }
\node at (3.5*\cs,\ylabel) {$R$};
\node at (3.5*\cs,8*\cs) {3};
\node at (4.5*\cs,8*\cs) {4};
\begin{scope}[shift={(8*\cs,0)}]
\fill[color0!30] (0,0) rectangle (7*\cs,7*\cs);
\foreach \x\y in {0/0, 1/1, 2/2, 3/3, 4/4, 5/5, 6/6, 5/3, 5/4}{
  \fill[color0] (\x*\cs,7*\cs-\y*\cs) rectangle (\x*\cs+\cs,6*\cs-\y*\cs);
}
\foreach \xy in {1,...,6}{
  \draw[gridline] (\xy*\cs,0)--(\xy*\cs,7*\cs);
  \draw[gridline] (0,\xy*\cs)--(7*\cs,\xy*\cs);
 }
\node at (3.5*\cs,\ylabel) {$V$};
\end{scope}
\end{scope}
\begin{scope}[shift={(11.1,0)}]
\fill[color0!30] (0,0) rectangle (7*\cs,7*\cs);
\fill[color2!30] (6*\cs,0) rectangle (7*\cs,7*\cs);
\fill[color2!30] (0,\cs) rectangle (7*\cs,2*\cs);
\foreach \x\y in {3/0, 3/2, 4/0, 4/1}{
  \fill[color0] (\x*\cs,7*\cs-\y*\cs) rectangle (\x*\cs+\cs,6*\cs-\y*\cs);
}
\foreach \x\y in {6/3, 6/4, 6/5}{
  \fill[color2] (\x*\cs,7*\cs-\y*\cs) rectangle (\x*\cs+\cs,6*\cs-\y*\cs);
}
\foreach \xy in {1,...,6}{
  \draw[gridline] (\xy*\cs,0)--(\xy*\cs,7*\cs);
  \draw[gridline] (0,\xy*\cs)--(7*\cs,\xy*\cs);
 }
\node at (3.5*\cs,\ylabel) {$R_c$};
\node at (6.5*\cs,8*\cs) {6};
\node at (-1*\cs,1.5*\cs) {5};
\begin{scope}[shift={(8*\cs,0)}]
\fill[color0!30] (0,0) rectangle (7*\cs,7*\cs);
\foreach \x\y in {0/0, 1/1, 2/2, 3/3, 4/4, 5/5, 6/6}{
  \fill[color0] (\x*\cs,7*\cs-\y*\cs) rectangle (\x*\cs+\cs,6*\cs-\y*\cs);
}
\foreach \xy in {1,...,6}{
  \draw[gridline] (\xy*\cs,0)--(\xy*\cs,7*\cs);
  \draw[gridline] (0,\xy*\cs)--(7*\cs,\xy*\cs);
 }
\node at (3.5*\cs,\ylabel) {$V_c$};
\end{scope}
\end{scope}
\end{tikzpicture}

%% file: figs/bigchanges.tex
\newcommand\rotdeg{18} 
\newcommand\leglength{1.3} 
\newcommand\matu{.1} 
\begin{tikzpicture}[
	celltype/.style={double=black,draw=white,double distance=1pt,line width=1.1pt, fill=color3, fill opacity=.2},
	edgetype/.style={double=color3,draw=white,double distance=1pt,line width=1.1pt},
  cellremove/.style={fill=color3, fill opacity=.2},
  cellkeep/.style={fill=black, fill opacity=.3},
  gridline/.style={white,line width=.3pt}
]
\foreach \x\label in {0/a, 3.8/b, 7.6/c, 11.1/d}{
  \node at (\x,-1.8) {(\label)};
}
\begin{scope}[shift={(0,0)}]
\coordinate (a) at (0,0);
\foreach \r\n in {0/1, 72/2, 144/3, 216/4, 288/5}{
  \coordinate (b\n) at ($(a)+(\r+\rotdeg:\leglength)$);
}
\draw[celltype] (b1)--(b2)--(b3)--(b4)--(b5)--(b1)--cycle;
\foreach \n in {b1, b2, b3, b4, b5}{
  \draw[edgetype] (a)--(\n);
  \fill[black] (\n) circle[radius=2pt];
}
\fill[color3] (a) circle[radius=2pt];
\end{scope}
\begin{scope}[shift={(3.8,0)}]
\coordinate (a) at (0,0);
\foreach \r\n in {0/1, 72/2, 144/3, 216/4, 288/5}{
  \coordinate (b\n) at ($(a)+(\r+\rotdeg:\leglength)$);
}
\fill[cellremove] (b4)--(a)--(b5)--cycle;
\fill[cellkeep] (b1) -- (b2) -- (b3) -- (b4) -- (a) -- (b5) -- (b1) -- (a) -- cycle;
\foreach \x\y in {b1/a, b1/b2, b2/a, b2/b3, b3/a, b4/a, b3/b4, b5/a, b5/b1}{
  \draw[celltype] (\x)--(\y);
} 
\draw[edgetype] (b4)--(b5);
\foreach \n in {a, b1, b2, b3, b4, b5}{
  \fill[black] (\n) circle[radius=2pt];
}
\end{scope}
\begin{scope}[shift={(7.6,-.1)}]
\coordinate (a) at (-1,1);
\coordinate (b) at (1,1);
\foreach \fac in {1,2,3,4,5}{
  \coordinate (c\fac) at ($(-1,-.5)+(0:\fac*.333)$);
  \draw[celltype] (a)--(c\fac)--(b);
}
\draw[edgetype] (a)--(b);
\foreach \n in {a, b, c1, c2, c3, c4, c5}{
  \fill[black] (\n) circle[radius=2pt];
}
\end{scope}
\begin{scope}[shift={(10.3,2)}]
\foreach \x\y\n in {.4/-.5/a, 0/-1.8/b, 1/-2/c, 1.2/-1.2/d}{
  \coordinate (\n) at (\x*1.5,\y*1.5);
}
\fill[cellkeep] (a)--(b)--(d);
\fill[cellkeep] (a)--(c)--(d);
\draw[celltype] (b)--(c)--(d);
\draw[celltype] (b)--(d);
\draw[celltype] (b)--(a)--(c);
\draw[celltype] (a)--(d);
\draw[edgetype] (b)--(c);
\foreach \n in {a, b, c,d}{
  \fill[black] (\n) circle[radius=2pt];
}
\end{scope}
\end{tikzpicture}

%% file: figs/reduced-matrices.tex
\newcommand\cwid{.4} 
\newcommand\subfigskip{3.3} 
\begin{tikzpicture}[gridline/.style={white,line width=.5pt}]
\fill[color0!30] (0,0) rectangle (4*\cwid,-4*\cwid);
\foreach \x\y in {0/0, 0/1, 1/1, 1/2, 2/2, 2/3}{
  \fill[color0] (\x*\cwid,-\y*\cwid) rectangle (\x*\cwid+\cwid,-\y*\cwid-\cwid);
}
\foreach \p\q in {1/{e_1}, 2/{e_2}, 3/{e_3}, 4/{e_4}}{
  \node[anchor=east,scale=.8] at (-.15,-\p*\cwid+.5*\cwid) {$v_{\p}$};
  \node[anchor=south,scale=.8] at (\p*\cwid-.5*\cwid,.15) {$\q$};
}
\foreach \xy in {1,2,3}{
  \draw[gridline] (\xy*\cwid,0)--(\xy*\cwid,-4*\cwid);
  \draw[gridline] (0,-1*\xy*\cwid)--(4*\cwid,-1*\xy*\cwid);
}
\node at (2*\cwid,-5*\cwid) {$R$};
\begin{scope}[shift={(-\subfigskip,-2.2)}]
\foreach \x\y\n in {1/1/1, 0/1/2, 0/0/3, 1/0/4}{ 
  \coordinate (\n) at (\x+.08,\y+1);
  \fill[black] (\n) circle[radius=2pt];
  \node[anchor={135+\n*90},inner sep=5pt] (v\n) at (\n) {$v_{\n}$};
}
\foreach \x\y in {1/2, 2/3, 3/4}{ 
  \draw[line width=1pt] (\x) -- (\y);
}
\node at ($(v1)!.5!(v2)$) {$e_1$};
\node at ($(v2)!.5!(v3)$) {$e_2$};
\node at ($(v3)!.5!(v4)$) {$e_3$};
\end{scope}
\begin{scope}[shift={(\subfigskip,0)}]
\node at (2*\cwid,2*\cwid) {$e_4 = \{v_1,v_3\}$};
\fill[color0!30] (0,0) rectangle (4*\cwid,-4*\cwid);
\foreach \x\y in {0/0, 1/1, 2/2, 3/3, 3/0, 3/1}{
  \fill[color0] (\x*\cwid,-\y*\cwid) rectangle (\x*\cwid+\cwid,-\y*\cwid-\cwid);
}
\node[anchor=east,scale=.8] at (-.15,-1*\cwid+.5*\cwid) {$e_1$};
\node[anchor=east,scale=.8] at (-.15,-2*\cwid+.5*\cwid) {$e_2$};
\node[anchor=south,scale=.8] at (4*\cwid-.5*\cwid,.15) {$e_4$};
\foreach \xy in {1,2,3}{
  \draw[gridline] (\xy*\cwid,0)--(\xy*\cwid,-4*\cwid);
  \draw[gridline] (0,-1*\xy*\cwid)--(4*\cwid,-1*\xy*\cwid);
}
\node at (2*\cwid,-5*\cwid) {$V_1$};
\end{scope}
\begin{scope}[shift={(2*\subfigskip,0)}]
\node at (2*\cwid,2*\cwid) {$e_4 = \{v_2,v_4\}$};
\fill[color0!30] (0,0) rectangle (4*\cwid,-4*\cwid);
\foreach \x\y in {0/0, 1/1, 2/2, 3/3, 3/1, 3/2}{
  \fill[color0] (\x*\cwid,-\y*\cwid) rectangle (\x*\cwid+\cwid,-\y*\cwid-\cwid);
}
\node[anchor=east,scale=.8] at (-.15,-2*\cwid+.5*\cwid) {$e_2$};
\node[anchor=east,scale=.8] at (-.15,-3*\cwid+.5*\cwid) {$e_3$};
\node[anchor=south,scale=.8] at (4*\cwid-.5*\cwid,.15) {$e_4$};
\foreach \xy in {1,2,3}{
  \draw[gridline] (\xy*\cwid,0)--(\xy*\cwid,-4*\cwid);
  \draw[gridline] (0,-1*\xy*\cwid)--(4*\cwid,-1*\xy*\cwid);
}
\node at (2*\cwid,-5*\cwid) {$V_2$};
\end{scope}
\begin{scope}[shift={(3*\subfigskip,0)}]
\node at (2*\cwid,2*\cwid) {$e_4 = \{v_1,v_4\}$};
\fill[color0!30] (0,0) rectangle (4*\cwid,-4*\cwid);
\foreach \x\y in {0/0, 1/1, 2/2, 3/3, 3/0, 3/1, 3/2}{
  \fill[color0] (\x*\cwid,-\y*\cwid) rectangle (\x*\cwid+\cwid,-\y*\cwid-\cwid);
}
\node[anchor=east,scale=.8] at (-.15,-1*\cwid+.5*\cwid) {$e_1$};
\node[anchor=east,scale=.8] at (-.15,-2*\cwid+.5*\cwid) {$e_2$};
\node[anchor=east,scale=.8] at (-.15,-3*\cwid+.5*\cwid) {$e_3$};
\node[anchor=south,scale=.8] at (4*\cwid-.5*\cwid,.15) {$e_4$};
\foreach \xy in {1,2,3}{
  \draw[gridline] (\xy*\cwid,0)--(\xy*\cwid,-4*\cwid);
  \draw[gridline] (0,-1*\xy*\cwid)--(4*\cwid,-1*\xy*\cwid);
}
\node at (2*\cwid,-5*\cwid) {$V_3$};
\end{scope}
\end{tikzpicture}

%% file: figs/diag_dey-hou-comparison0.tex
\newcommand\mainrad{.8} 
\setlength{\tabcolsep}{0pt}
\raisebox{1.45cm}{$
\begin{array}{c C c C c C c C c C c C c C c C c}
&&
\begin{tikzpicture}[baseline=0pt]\useasboundingbox (-.1,0) rectangle (.1,3.5*\mainrad);
\fill[black] (0,0) circle[radius=2pt];
\node[anchor=south,inner sep=6pt] at (0,0) {$v_1$};
\end{tikzpicture}
&&
\begin{tikzpicture}[baseline=0pt]\useasboundingbox (-.1,0) rectangle (.1,3.5*\mainrad);
\fill[black] (0,0) circle[radius=2pt];
\fill[black] (0,\mainrad) circle[radius=2pt];
\node[anchor=south,inner sep=6pt] at (0,\mainrad) {$v_2$};
\end{tikzpicture}
&&
\begin{tikzpicture}[baseline=0pt]\useasboundingbox (-.1,0) rectangle (.1,3.5*\mainrad);
\fill[black] (0,0) circle[radius=2pt];
\fill[black] (0,\mainrad) circle[radius=2pt];
\fill[black] (0,2*\mainrad) circle[radius=2pt];
\node[anchor=south,inner sep=6pt] at (0,2*\mainrad) {$v_3$};
\end{tikzpicture}
&&
\begin{tikzpicture}[baseline=0pt]\useasboundingbox (-.1,0) rectangle (.1,3.5*\mainrad);
\fill[black] (0,0) circle[radius=2pt];
\fill[black] (0,\mainrad) circle[radius=2pt];
\fill[black] (0,2*\mainrad) circle[radius=2pt];
\fill[black] (0,3*\mainrad) circle[radius=2pt];
\node[anchor=south,inner sep=6pt] at (0,3*\mainrad) {$v_4$};
\end{tikzpicture}
&&
\begin{tikzpicture}[baseline=0pt]\useasboundingbox (-.1,0) rectangle (.1,3.5*\mainrad);
\draw[line width=1pt,color3] (0,0) to node[left,black] {$e_1$} node[right,white] {$e_1$} (0,\mainrad);
\fill[black] (0,0) circle[radius=2pt];
\fill[black] (0,\mainrad) circle[radius=2pt];
\fill[black] (0,2*\mainrad) circle[radius=2pt];
\fill[black] (0,3*\mainrad) circle[radius=2pt];
\end{tikzpicture}
&&
\begin{tikzpicture}[baseline=0pt]\useasboundingbox (-.1,0) rectangle (.1,3.5*\mainrad);
\draw[line width=1pt,color3] (0,0) -- (0,\mainrad);
\draw[line width=1pt] (0,\mainrad) to node[left] {$e_2$} node[right,white] {$e_2$} (0,2*\mainrad);
\fill[black] (0,0) circle[radius=2pt];
\fill[black] (0,\mainrad) circle[radius=2pt];
\fill[black] (0,2*\mainrad) circle[radius=2pt];
\fill[black] (0,3*\mainrad) circle[radius=2pt];
\end{tikzpicture}
&&
\begin{tikzpicture}[baseline=0pt]\useasboundingbox (-.1,0) rectangle (.1,3.5*\mainrad);
\draw[line width=1pt,color3] (0,0) -- (0,\mainrad);
\draw[line width=1pt] (0,\mainrad) to node[left,pos=.75] {$e_3$} node[right,white] {$e_3$} (0,3*\mainrad);
\fill[black] (0,0) circle[radius=2pt];
\fill[black] (0,\mainrad) circle[radius=2pt];
\fill[black] (0,2*\mainrad) circle[radius=2pt];
\fill[black] (0,3*\mainrad) circle[radius=2pt];
\end{tikzpicture}
\\
\emptyset & \xhookrightarrow{v_1} &
K_1 & \hspace{-.1cm} \xhookrightarrow{v_2} &
K_2 & \hspace{-.1cm} \xhookrightarrow{v_3} & 
K_3 & \hspace{-.1cm} \xhookrightarrow{v_4} & 
K_4 & \hspace{-.1cm} \xhookrightarrow{e_1} & 
K_5 & \hspace{-.1cm} \xhookrightarrow{e_2} & 
K_6 & \hspace{-.1cm} \xhookrightarrow{e_3} & 
K_7 & \hspace{-.1cm} \xhookleftarrow{e_3} & \cdots.
\end{array}
$}

%% file: figs/diag_dey-hou-comparison1.tex
\newcommand\mainrad{.8} 
\newcommand\secrad{.5} 
\newcommand\bobox{.35} 
\setlength{\tabcolsep}{0pt}
\raisebox{1.33cm}{$
\begin{array}{c C c C c C c C c C c C c C c C c}
&&
\begin{tikzpicture}[baseline=0pt]\useasboundingbox (-1*\bobox,0) rectangle (\bobox,3*\mainrad);
\draw[line width=1pt,color2] (0,0) -- (0,\mainrad);
\fill[black] (0,0) circle[radius=2pt];
\fill[black] (0,\mainrad) circle[radius=2pt];
\fill[black] (0,2*\mainrad) circle[radius=2pt];
\fill[black] (0,3*\mainrad) circle[radius=2pt];
\fill[black] (\secrad,1.5*\mainrad) circle[radius=2pt];
\node[anchor=south,inner sep=6pt] at (\secrad,1.5*\mainrad) {$\omega$};
\draw[line width=1pt] (0,\mainrad)--(0,3*\mainrad);
\end{tikzpicture}
&&
\begin{tikzpicture}[baseline=0pt]\useasboundingbox (-1*\bobox,0) rectangle (\bobox,3*\mainrad);
\draw[line width=1pt,color2] (0,0) -- (0,\mainrad);
\fill[black] (0,0) circle[radius=2pt];
\fill[black] (0,\mainrad) circle[radius=2pt];
\fill[black] (0,2*\mainrad) circle[radius=2pt];
\fill[black] (0,3*\mainrad) circle[radius=2pt];
\fill[black] (\secrad,1.5*\mainrad) circle[radius=2pt];
\draw[line width=1pt] (0,\mainrad)--(0,3*\mainrad);
\draw[line width=1pt] (0,0)--(\secrad,1.5*\mainrad);
\end{tikzpicture}
&&
\begin{tikzpicture}[baseline=0pt]\useasboundingbox (-1*\bobox,0) rectangle (\bobox,3*\mainrad);
\draw[line width=1pt,color2] (0,0) -- (0,\mainrad);
\fill[black] (0,0) circle[radius=2pt];
\fill[black] (0,\mainrad) circle[radius=2pt];
\fill[black] (0,2*\mainrad) circle[radius=2pt];
\fill[black] (0,3*\mainrad) circle[radius=2pt];
\fill[black] (\secrad,1.5*\mainrad) circle[radius=2pt];
\draw[line width=1pt] (0,\mainrad)--(0,3*\mainrad);
\draw[line width=1pt] (0,0)--(\secrad,1.5*\mainrad);
\draw[line width=1pt] (0,\mainrad)--(\secrad,1.5*\mainrad);
\end{tikzpicture}
&&
\begin{tikzpicture}[baseline=0pt]\useasboundingbox (-1*\bobox,0) rectangle (\bobox,3*\mainrad);
\draw[line width=1pt,color2] (0,0) -- (0,\mainrad);
\fill[black] (0,0) circle[radius=2pt];
\fill[black] (0,\mainrad) circle[radius=2pt];
\fill[black] (0,2*\mainrad) circle[radius=2pt];
\fill[black] (0,3*\mainrad) circle[radius=2pt];
\fill[black] (\secrad,1.5*\mainrad) circle[radius=2pt];
\draw[line width=1pt] (0,\mainrad)--(0,3*\mainrad);
\draw[line width=1pt] (0,0)--(\secrad,1.5*\mainrad);
\draw[line width=1pt] (0,\mainrad)--(\secrad,1.5*\mainrad);
\draw[line width=1pt] (0,2*\mainrad)--(\secrad,1.5*\mainrad);
\end{tikzpicture}
&&
\begin{tikzpicture}[baseline=0pt]\useasboundingbox (-1*\bobox,0) rectangle (\bobox,3*\mainrad);
\draw[line width=1pt,color2] (0,0) -- (0,\mainrad);
\fill[black] (0,0) circle[radius=2pt];
\fill[black] (0,\mainrad) circle[radius=2pt];
\fill[black] (0,2*\mainrad) circle[radius=2pt];
\fill[black] (0,3*\mainrad) circle[radius=2pt];
\fill[black] (\secrad,1.5*\mainrad) circle[radius=2pt];
\draw[line width=1pt] (0,\mainrad)--(0,3*\mainrad);
\draw[line width=1pt] (0,0)--(\secrad,1.5*\mainrad);
\draw[line width=1pt] (0,\mainrad)--(\secrad,1.5*\mainrad);
\draw[line width=1pt] (0,2*\mainrad)--(\secrad,1.5*\mainrad);
\draw[line width=1pt] (0,3*\mainrad)--(\secrad,1.5*\mainrad);
\end{tikzpicture}
&&
\begin{tikzpicture}[baseline=0pt]\useasboundingbox (-1*\bobox,0) rectangle (\bobox,3*\mainrad);
\fill[color2!20] (0,0) -- (0,\mainrad) -- (\secrad,1.5*\mainrad) -- cycle;
\draw[line width=1pt,color2] (0,0) -- (0,\mainrad);
\fill[black] (0,0) circle[radius=2pt];
\fill[black] (0,\mainrad) circle[radius=2pt];
\fill[black] (0,2*\mainrad) circle[radius=2pt];
\fill[black] (0,3*\mainrad) circle[radius=2pt];
\fill[black] (\secrad,1.5*\mainrad) circle[radius=2pt];
\draw[line width=1pt] (0,\mainrad)--(0,3*\mainrad);
\draw[line width=1pt] (0,0)--(\secrad,1.5*\mainrad);
\draw[line width=1pt] (0,\mainrad)--(\secrad,1.5*\mainrad);
\draw[line width=1pt] (0,2*\mainrad)--(\secrad,1.5*\mainrad);
\draw[line width=1pt] (0,3*\mainrad)--(\secrad,1.5*\mainrad);
\end{tikzpicture}
&&
\begin{tikzpicture}[baseline=0pt]\useasboundingbox (-1*\bobox,0) rectangle (\bobox,3*\mainrad);
\fill[color2!20] (0,0) -- (0,\mainrad) -- (\secrad,1.5*\mainrad) -- cycle;
\fill[black!20] (0,\mainrad) -- (0,2*\mainrad) -- (\secrad,1.5*\mainrad) -- cycle;
\draw[line width=1pt,color2] (0,0) -- (0,\mainrad);
\fill[black] (0,0) circle[radius=2pt];
\fill[black] (0,\mainrad) circle[radius=2pt];
\fill[black] (0,2*\mainrad) circle[radius=2pt];
\fill[black] (0,3*\mainrad) circle[radius=2pt];
\fill[black] (\secrad,1.5*\mainrad) circle[radius=2pt];
\draw[line width=1pt] (0,\mainrad)--(0,3*\mainrad);
\draw[line width=1pt] (0,0)--(\secrad,1.5*\mainrad);
\draw[line width=1pt] (0,\mainrad)--(\secrad,1.5*\mainrad);
\draw[line width=1pt] (0,2*\mainrad)--(\secrad,1.5*\mainrad);
\draw[line width=1pt] (0,3*\mainrad)--(\secrad,1.5*\mainrad);
\end{tikzpicture}
&&
\begin{tikzpicture}[baseline=0pt]\useasboundingbox (-1*\bobox,0) rectangle (\bobox,3*\mainrad);
\fill[color2!20] (0,0) -- (0,\mainrad) -- (\secrad,1.5*\mainrad) -- cycle;
\fill[black!20] (0,\mainrad) -- (0,3*\mainrad) -- (\secrad,1.5*\mainrad) -- cycle;
\draw[line width=1pt,color2] (0,0) -- (0,\mainrad);
\fill[black] (0,0) circle[radius=2pt];
\fill[black] (0,\mainrad) circle[radius=2pt];
\fill[black] (0,2*\mainrad) circle[radius=2pt];
\fill[black] (0,3*\mainrad) circle[radius=2pt];
\fill[black] (\secrad,1.5*\mainrad) circle[radius=2pt];
\draw[line width=1pt] (0,\mainrad)--(0,3*\mainrad);
\draw[line width=1pt] (0,0)--(\secrad,1.5*\mainrad);
\draw[line width=1pt] (0,\mainrad)--(\secrad,1.5*\mainrad);
\draw[line width=1pt] (0,2*\mainrad)--(\secrad,1.5*\mainrad);
\draw[line width=1pt] (0,3*\mainrad)--(\secrad,1.5*\mainrad);
\end{tikzpicture}
\\
\cdots & \xhookrightarrow{e_3} &
\widehat{K}_7 & \hspace{-.1cm} \xhookrightarrow{\widehat{v_1}} &
\widehat{K}_8 & \hspace{-.1cm} \xhookrightarrow{\widehat{v_2}} &
\widehat{K}_9 & \hspace{-.1cm} \xhookrightarrow{\widehat{v_3}} &
\widehat{K}_{10} & \hspace{-.1cm} \xhookrightarrow{\widehat{v_4}} &
\widehat{K}_{11} & \hspace{-.1cm} \xhookrightarrow{\widehat{e_1}} &
\widehat{K}_{12} & \hspace{-.1cm} \xhookrightarrow{\widehat{e_2}} &
\widehat{K}_{13} & \hspace{-.1cm} \xhookrightarrow{\widehat{e_3}} &
\widehat{K}_{14}.
\end{array}
$}

%% file: figs/diag_dey-hou-comparison2.tex
\newcommand\mainrad{1} 
\newcommand\cwid{.2} 
\newcommand\mainskip{45}
\newcommand\graycol{black} 
D\ = \ \begin{tikzpicture}[
  baseline=-4*\cwid cm,
  gridline/.style={white,line width=.5pt}
]
\fill[color0!30] (0,0) rectangle (7*\cwid,-7*\cwid);
\fill[color3!30] (4*\cwid,0) rectangle (5*\cwid,-7*\cwid);
\fill[color3!30] (0,-4*\cwid) rectangle (7*\cwid,-5*\cwid);
\fill[color3] (4*\cwid,0) rectangle (5*\cwid,-2*\cwid);
\fill[color0] (5*\cwid,-1*\cwid) rectangle (6*\cwid,-3*\cwid);
\fill[color0] (6*\cwid,-2*\cwid) rectangle (7*\cwid,-4*\cwid);
\foreach \xy in {1,...,6}{
  \draw[gridline] (\xy*\cwid,0)--(\xy*\cwid,-7*\cwid);
  \draw[gridline] (0,-1*\xy*\cwid)--(7*\cwid,-1*\xy*\cwid);
}
\node at (4.5*\cwid, 1.5*\cwid) {$e_1$};
\end{tikzpicture}
\ \ , \ 
V\ =\ 
\begin{tikzpicture}[
  baseline=-4*\cwid cm,
  gridline/.style={white,line width=.5pt}
]
\fill[color0!30] (0,0) rectangle (7*\cwid,-7*\cwid);
\fill[color3!30] (4*\cwid,0) rectangle (5*\cwid,-7*\cwid);
\fill[color3!30] (0,-4*\cwid) rectangle (7*\cwid,-5*\cwid);
\foreach \xy in {0,1,2,3,5,6}{
  \fill[color0] (\xy*\cwid,-1*\xy*\cwid) rectangle (\xy*\cwid+\cwid,-1*\xy*\cwid-\cwid);
}
\fill[color3] (4*\cwid,-4*\cwid) rectangle (5*\cwid,-5*\cwid);
\foreach \xy in {1,...,6}{
  \draw[gridline] (\xy*\cwid,0)--(\xy*\cwid,-7*\cwid);
  \draw[gridline] (0,-1*\xy*\cwid)--(7*\cwid,-1*\xy*\cwid);
}
\end{tikzpicture}
\ \ , \ 
\widehat D\ =\ 
\begin{tikzpicture}[
  baseline=-8*\cwid cm,
  gridline/.style={white,line width=.5pt}
]
\fill[color0!30] (0,0) rectangle (15*\cwid,-15*\cwid);
\fill[color2!30] (5*\cwid,0) rectangle (6*\cwid,-15*\cwid);
\fill[color2!30] (0,-5*\cwid,0) rectangle (15*\cwid,-6*\cwid);
\fill[color2!30] (12*\cwid,0) rectangle (13*\cwid,-15*\cwid);
\fill[color2!30] (0,-12*\cwid) rectangle (15*\cwid,-13*\cwid);
\foreach \x\y in {
7/3, 7/4, 8/4, 8/5, 9/1, 9/2, 10/1, 10/3, 11/1, 11/4, 12/1, 12/5,
14/7, 14/10, 14/11, 15/8, 15/11, 15/12
}{
  \fill[color0] (\x*\cwid-\cwid,-1*\y*\cwid+\cwid) rectangle (\x*\cwid,-1*\y*\cwid);
}
\foreach \x\y in {6/2, 6/3, 13/6, 13/9, 13/10}{
  \fill[color2] (\x*\cwid-\cwid,-1*\y*\cwid+\cwid) rectangle (\x*\cwid,-1*\y*\cwid);
}
\foreach \xy in {1,...,14}{
  \draw[gridline] (\xy*\cwid,0)--(\xy*\cwid,-15*\cwid);
  \draw[gridline] (0,-1*\xy*\cwid)--(15*\cwid,-1*\xy*\cwid);
}
\node at (5.5*\cwid, 1.5*\cwid) {$e_1$};
\node at (12.5*\cwid, 1.5*\cwid) {$\widehat{e_1}$};
\end{tikzpicture}
\ \ ,\ 
\widehat V\ =\ 
\begin{tikzpicture}[
  baseline=-8*\cwid cm,
  gridline/.style={white,line width=.5pt}
]
\fill[color0!30] (0,0) rectangle (15*\cwid,-15*\cwid);
\fill[color2!30] (5*\cwid,0) rectangle (6*\cwid,-15*\cwid);
\fill[color2!30] (0,-5*\cwid,0) rectangle (15*\cwid,-6*\cwid);
\fill[color2!30] (12*\cwid,0) rectangle (13*\cwid,-15*\cwid);
\fill[color2!30] (0,-12*\cwid) rectangle (15*\cwid,-13*\cwid);
\foreach \x\y in {
1/1, 2/2, 3/3, 4/4, 5/5, 7/7, 8/8, 9/9, 10/10, 11/11, 12/12, 14/14, 15/15,
10/9, 11/7, 11/9, 12/8, 12/7, 12/9
}{
  \fill[color0] (\x*\cwid-\cwid,-1*\y*\cwid+\cwid) rectangle (\x*\cwid,-1*\y*\cwid);
}
\foreach \x\y in {
6/6, 13/13,
10/6, 11/6, 12/6
}{
  \fill[color2] (\x*\cwid-\cwid,-1*\y*\cwid+\cwid) rectangle (\x*\cwid,-1*\y*\cwid);
}
\foreach \xy in {1,...,14}{
  \draw[gridline] (\xy*\cwid,0)--(\xy*\cwid,-15*\cwid);
  \draw[gridline] (0,-1*\xy*\cwid)--(15*\cwid,-1*\xy*\cwid);
}
\end{tikzpicture}

%% file: figs/testing_environment.tex
\newcommand\pich{4.2cm} 
\newcommand\picsh{.1} 
\newcommand\picshh{.15} 
\newcommand\colwid{3.3} 
\newcommand\colwidtwo{9.8} 
\newcommand\outputoffset{.5} 
\begin{tikzpicture}[
	flowbig/.style={rectangle,draw,minimum width=4cm, minimum height=.7cm},
	flowsmall/.style={rectangle,draw,minimum width=2cm, minimum height=.7cm},
	outof/.style={line width=1pt,-{Straight Barb[width=4pt,length=2pt]},shorten <=4pt},
	outoftrunc/.style={line width=1pt,-{Straight Barb[width=4pt,length=2pt]},shorten >=2pt},
	into/.style={line width=1pt,{Bar[width=4pt]}-{Straight Barb[width=4pt,length=2pt]}, shorten >=2pt},
	intotrunc/.style={line width=1pt,{Bar[width=4pt]}-},
	conx/.style={line width=1pt}
]
\fill[color7!15] (-.8,2.8) rectangle (8,.5);
\fill[color7!15] (8.5,2.8) rectangle (14.35,.5);
\node[anchor=south west] at (-.75,2.8) {\tiny\texttt{Preparation for testing}};
\node[anchor=south west] at (8.5,2.8) {\tiny\texttt{Testing}};
\node[flowsmall,fill=color2!30] (gen) at (0.4,2.3) {\tiny\texttt{generate.ipynb\vphantom{Ap}}};
\node (gen_output) at (.4,1) {\tiny\texttt{dataset.bmat\vphantom{Ap}}};
\node[flowsmall,fill=color2!30] (update) at (\colwid,2.3) {\tiny\texttt{choose simplices\vphantom{Ap}}};
\node[anchor=west] (update_output) at ($(update.east)+(0:\outputoffset)$) {\tiny\parbox{2.65cm}{\raggedright\texttt{dataset\_updated.bmat}\\\texttt{dataset.update}}};
\node[flowsmall,fill=color0!30] (phatop1) at (\colwid,1) {\tiny\texttt{PHAT-op\vphantom{Ap}}};
\node[anchor=west] (phatop1_output) at ($(phatop1.east)+(0:\outputoffset)$) {\tiny\parbox{2.65cm}{\raggedright\texttt{dataset\_reduced.bmat}\\\texttt{dataset.opmat}}};
\node[flowsmall,fill=color0!30] (phatop2) at (\colwidtwo,2.3) {\tiny\texttt{PHAT-op\vphantom{Ap}}};
\node[anchor=west] (phatop2_output) at ($(phatop2.east)+(0:\outputoffset)$) {\tiny\parbox{2.3cm}{\raggedright\texttt{dataset\_reduced.bmat}\\\texttt{dataset.opmat}}};
\node[flowsmall,fill=color1!30] (sirup) at (\colwidtwo,1) {\tiny\texttt{SiRUP\vphantom{Ap}}};
\node[anchor=west] (sirup_output) at ($(sirup.east)+(0:\outputoffset)$) {\tiny\parbox{2.3cm}{\raggedright\texttt{dataset\_reduced.bmat}\\\texttt{dataset.opmat}}};
\draw[outof] (gen) -- (gen_output);
\draw[outof,transform canvas={yshift=.13cm}] (phatop2) -- (phatop2_output);
\draw[outof,transform canvas={yshift=-.11cm}] (phatop2) -- (phatop2_output);
\draw[outof,transform canvas={yshift=.13cm}] (sirup) -- (sirup_output);
\draw[outof,transform canvas={yshift=-.11cm}] (sirup) -- (sirup_output);
\draw[outof,transform canvas={yshift=.13cm}] (update) -- (update_output);
\draw[outof,transform canvas={yshift=-.11cm}] (update) -- (update_output);
\draw[outof,transform canvas={yshift=.13cm}] (phatop1) -- (phatop1_output);
\draw[outof,transform canvas={yshift=-.11cm}] (phatop1) -- (phatop1_output);
\draw[outoftrunc] ($(gen_output.east)!.5!(phatop1.west)$) |- (update.west);
\draw[into] (gen_output.east) -- (phatop1.west);
\draw[into,transform canvas={yshift=.13cm}] ($(update_output.east)+(180:.1)$) -- (phatop2.west);
\draw[intotrunc,transform canvas={yshift=-.11cm}] ($(update_output.east)+(180:.75)$) -| ($($(phatop1_output.east)+(90:.24)$)!.5!($(sirup.west)+(90:.24)$)$);
\draw[into,transform canvas={yshift=.13cm}] ($(phatop1_output.east)+(180:.1)$) --  (sirup.west);
\draw[into,transform canvas={yshift=-.11cm}] ($(phatop1_output.east)+(180:.9)$) --  (sirup.west);
\end{tikzpicture}

%% file: Bibliography.bib
@article{keepitsparse, 
	title={Keeping it sparse: Computing Persistent Homology revisited}, 
	volume={3}, 
	number={1}, 
	journal={Computing in Geometry and Topology}, 
	author={Bauer, Ulrich and Bin Masood, Talha and Giunti, Barbara and Houry, Guillaume and Kerber, Michael and Rathod , Abhishek}, 
	year={2024}, 
	pages={6:1–6:26} 
}

@ARTICLE{cycles_userguide,
	title    = {Minimal Cycle Representatives in Persistent Homology Using Linear
	Programming: An Empirical Study With User's Guide},
	author   = {Li, Lu and Thompson, Connor and Henselman-Petrusek, Gregory and
	Giusti, Chad and Ziegelmeier, Lori},
	journal  = {Front Artif Intell},
	volume   = {4},
	pages    = {681117},
	year     = {2021},
}

@article{tight_cycles,
	title    = {Tight basis cycle representatives for persistent homology of large biological data sets},
	author   = {Aggarwal, Manu and Periwal, Vipul},
	journal  = {PLoS Comput Biol},
	volume   =  {19},
	number   =  {5},
	pages    = "{e1010341}",
	year     =  {2023},
}

@misc{hickok_bundle,
	title={Persistence Diagram Bundles: A Multidimensional Generalization of Vineyards}, 
	author={Abigail Hickok},
	year={2023},
	note={Preprint available at {arXiv}:2210.05124},
}

@InProceedings{cerrifrosini,
	author={Cerri, Andrea and Ethier, Marc and Frosini, Patrizio},
	title={A Study of Monodromy in the Computation of Multidimensional Persistence},
	booktitle={Discrete Geometry for Computer Imagery},
	year={2013},
	publisher={Springer Berlin Heidelberg},
	pages={192--202},
}

@INPROCEEDINGS{notes_pivot,
	author = {Giunti, Barbara},
	title = {Notes on pivot pairings},
	booktitle = {Proceedings of the 37th European Workshop on Computational Geometry},
	year = {2021},
	pages = {11:1--11:6}
}

@inproceedings{adams2014javaplex,
 	title={JavaPlex: A research software package for persistent (co)homology},	
	 author={Adams, Henry and Tausz, Andrew and Vejdemo-Johansson, Mikael},
 	booktitle={Mathematical Software -- ICMS 2014},
 	year={2014},
 	publisher={Springer Berlin Heidelberg},
 	address={Lecture Notes in Computer Science},
 	volume={8592},
 	pages={129--136},
}

@misc{dey2022updating,
	title={Updating Barcodes and Representatives for Zigzag Persistence}, 
	author={Tamal K. Dey and Tao Hou},
	year={2022},
	note={Preprint available at arXiv:2112.02352},
}

@misc{donut,
	author = {Giunti, Barbara and Lazovskis, J{\=a}nis and Rieck, Bastian},
	title  = {
	{DONUT}: {D}atabase of {O}riginal \& {N}on-{T}heoretical {U}ses of {T}opology},
	note   = {\url{https://donut.topology.rocks}},
	year   = {2022},
	key    = {DONUT},
}

@misc{tadasets,
	title = {{TaDAset - a Scikit-TDA project}},
  author = {Nathaniel Saul and Chris Tralie},
	year = {2018},
	publisher = {GitHub},
	journal = {GitHub repository},
	howpublished = {\url{https://github.com/scikit-tda/tadasets}},
	key = {Tadaset},
	doi = {10.5281/zenodo.2533369},
}

@misc{our_zenodo,
  author = {Giunti, Barbara and Lazovskis, Jānis},
  title = {Dataset for ``{P}runing vineyards: updating barcodes and representative cycles by removing simplices"},
  year = {2025},
  publisher = {Zenodo},
  note = {\url{https://zenodo.org/records/15481043}},
}

@misc{images_database,
	author = {Klacansky, Pavol},
	title = {{Open Scientific Visualization Datasets}},
	howpublished = {\url{https://klacansky.com/open-scivis-datasets/}},
	key = {Database images},
}

@InProceedings{fast_zigzag,
  author =	{Dey, Tamal K. and Hou, Tao},
  title =	{{Fast Computation of Zigzag Persistence}},
  booktitle =	{30th Annual European Symposium on Algorithms (ESA 2022)},
  pages =	{43:1--43:15},
  series =	{Leibniz International Proceedings in Informatics (LIPIcs)},
  year =	{2022},
  volume =	{244},
  publisher =	{Schloss Dagstuhl -- Leibniz-Zentrum f{\"u}r Informatik},
}

@InProceedings{dey2023computing,
	author =	{Dey, Tamal K. and Hou, Tao},
	title =	{{Computing Zigzag Vineyard Efficiently Including Expansions and Contractions}},
	booktitle =	{40th International Symposium on Computational Geometry (SoCG 2024)},
	pages =	{49:1--49:15},
	series =	{Leibniz International Proceedings in Informatics (LIPIcs)},
	year =	{2024},
	publisher =	{Schloss Dagstuhl -- Leibniz-Zentrum f{\"u}r Informatik},
}

@misc{ripser_software,
	author = {Ulrich Bauer},
	title = {Ripser},
	year = {2021},
	howpublished = {\url{https://github.com/Ripser/ripser}},
	note = {Version 1.2.1}
}

@inproceedings{standard_alg,
	title={Topological persistence and simplification},
	author={Edelsbrunner, Herbert and Letscher, David and Zomorodian, Afra},
	booktitle={Proceedings 41st annual symposium on foundations of computer science},
	pages={454--463},
	year={2000},
	organization={IEEE},
}

@misc{phat_vineyards,
	author = {Barbara Giunti and J\={a}nis Lazovskis},
	title = {\textsc{PHAT-vineyards}},
	year = {2024},
	publisher = {Bitbucket},
	journal = {Bitbucket repository},
	howpublished = {\url{https://bitbucket.org/jlazovskis/phat-sirup}}
}

@article{luo2023accelerating,
title = {Accelerating iterated persistent homology computations with warm starts},
journal = {Computational Geometry},
volume = {120},
pages = {102089},
year = {2024},
author = {Yuan Luo and Bradley J. Nelson},
}

@article{roadmap,
	title={A roadmap for the computation of persistent homology},
	author={Otter, Nina and Porter, Mason A and Tillmann, Ulrike and Grindrod, Peter and Harrington, Heather A},
	journal={EPJ Data Science},
	volume={6},
	pages={1--38},
	year={2017},
	publisher={Springer},
}

@book {oudot,
  AUTHOR = {Oudot, Steve Y.},
  TITLE = {Persistence theory: from quiver representations to data analysis},
	SERIES = {Mathematical Surveys and Monographs},
	VOLUME = {209},
	PUBLISHER = {American Mathematical Society, Providence, RI},
  YEAR = {2015},
 }

@book{edels_harer,
	AUTHOR = {Edelsbrunner, Herbert and Harer, John L.},
	TITLE = {Computational topology: An introduction},
	PUBLISHER = {American Mathematical Society, Providence, RI},
	YEAR = {2010},
	PAGES = {xii+241},
	volume={69}
}

@book{dey_wang_2022,
	place={Cambridge}, 
	title={Computational Topology for Data Analysis}, 
	publisher={Cambridge University Press}, 
	author={Dey, Tamal Krishna and Wang, Yusu}, 
	year={2022},
}

@inproceedings{vines_and_vineyards,
	author = {Cohen-Steiner, David and Edelsbrunner, Herbert and Morozov, Dmitriy},
	title = {Vines and Vineyards by Updating Persistence in Linear Time},
	year = {2006},
	publisher = {Association for Computing Machinery},
	address = {New York, NY, USA},
	booktitle = {Proceedings of the Twenty-Second Annual Symposium on Computational Geometry},
	pages = {119–126},
}

@InProceedings{facundo_dynamic,
author={Kim, Woojin and M{\'e}moli, Facundo and Smith, Zane},
title={Analysis of Dynamic Graphs and Dynamic Metric Spaces via Zigzag Persistence},
booktitle={Topological Data Analysis},
year={2020},
publisher={Springer International Publishing},
pages={371--389},
}

@incollection{clearcompress,
	author    = {Ulrich Bauer and Michael Kerber and Jan Reininghaus},
	title     = {Clear and Compress: Computing Persistent Homology in Chunks},
	booktitle = {Topological Methods in Data Analysis and Visualization III, Theory, Algorithms, and Applications},
	pages     = {103--117},
	publisher = {Springer},
	year      = {2014},
}

@INPROCEEDINGS{twist,
	author = {Chao Chen and Michael Kerber},
	title = {Persistent homology computation with a twist},
	booktitle = {Proceedings 27th European Workshop on Computational Geometry},
	year = {2011}
}

@misc{dionysus,
	title={Dionysus, a {C}++ library for computing persistent homology},
	author={Morozov, Dmitriy},
	year={2007}
}

@InProceedings{matrix_multi,
  author =	{Morozov, Dmitriy and Skraba, Primoz},
  title =	{{Persistent (Co)Homology in Matrix Multiplication Time}},
  booktitle =	{41st International Symposium on Computational Geometry (SoCG 2025)},
  pages =	{68:1--68:16},
  series =	{Leibniz International Proceedings in Informatics (LIPIcs)},
  year =	{2025},
  volume =	{332},
  editor =	{Aichholzer, Oswin and Wang, Haitao},
  publisher =	{Schloss Dagstuhl -- Leibniz-Zentrum f{\"u}r Informatik},
  address =	{Dagstuhl, Germany},
  doi =		{10.4230/LIPIcs.SoCG.2025.68},
}

@article{dualities_persistent,
	year = {2011},
	publisher = {{IOP} Publishing},
	pages = {124003},
	author = {Vin de Silva and Dmitriy Morozov and Mikael Vejdemo-Johansson},
	title = {Dualities in persistent (co)homology},
	journal = {Inverse Problems},
}

@article{bauer2021ripser,
	title={Ripser: efficient computation of {V}ietoris--{R}ips persistence barcodes},
	author={Bauer, Ulrich},
	journal={Journal of Applied and Computational Topology},
	pages={1--33},
	year={2021},
	volume={5},
 	number={3},
	publisher={Springer}
}

@article{phat_paper,
	title={Phat--persistent homology algorithms toolbox},
	author={Bauer, Ulrich and Kerber, Michael and Reininghaus, Jan and Wagner, Hubert},
	journal={Journal of Symbolic Computation},
 	volume={78},
  pages={76--90},
  year={2017},
  publisher={Elsevier}
}

@misc{henselman2016matroid,
 	title={Matroid filtrations and computational persistent homology},
 	author={Henselman, Gregory and Ghrist, Robert},
 	note={Preprint available at {arXiv}:1606.00199},
 	year={2016}
}

@incollection{maria2014gudhi
, author    = {Cl{\'{e}}ment Maria}
, title     = {Filtered Complexes}
, publisher = {GUDHI Editorial Board}
, edition   = {3.11.0}
, booktitle = {GUDHI User and Reference Manual}
, url       = {https://gudhi.inria.fr/doc/3.11.0/group__simplex__tree.html}
, year      = {2025}
}

@misc{perez2021giotto,
	title={Giotto-ph: A Python Library for High-Performance Computation of Persistent Homology of {V}ietoris--{R}ips Filtrations},
	author={P{\'e}rez, Juli{\'a}n Burella and Hauke, Sydney and Lupo, Umberto and Caorsi, Matteo and Dassatti, Alberto},
	note={Preprint available at {arXiv}:2107.05412},
	year={2021}
}

@ARTICLE{ccliques,
	AUTHOR={Reimann, Michael W.  and Nolte, Max  and Scolamiero, Martina  and Turner, Katharine  and Perin, Rodrigo  and Chindemi, Giuseppe  and Dłotko, Paweł  and Levi, Ran  and Hess, Kathryn  and Markram, Henry },
	TITLE={Cliques of Neurons Bound into Cavities Provide a Missing Link between Structure and Function},
	JOURNAL={Frontiers in Computational Neuroscience},
	VOLUME={11},
	YEAR={2017},
	ISSN={1662-5188},
}

@article{bb,
title = {Reconstruction and Simulation of Neocortical Microcircuitry},
journal = {Cell},
volume = {163},
number = {2},
pages = {456-492},
year = {2015},
author= {Markram, Henry and Muller, Eilif and Ramaswamy, Srikanth and
Reimann, Michael W and Abdellah, Marwan and Sanchez, Carlos
Aguado and Ailamaki, Anastasia and Alonso-Nanclares, Lidia and
Antille, Nicolas and Arsever, Selim and Kahou, Guy Antoine
Atenekeng and Berger, Thomas K and Bilgili, Ahmet and Buncic,
Nenad and Chalimourda, Athanassia and Chindemi, Giuseppe and
Courcol, Jean-Denis and Delalondre, Fabien and Delattre, Vincent
and Druckmann, Shaul and Dumusc, Raphael and Dynes, James and
Eilemann, Stefan and Gal, Eyal and Gevaert, Michael Emiel and
Ghobril, Jean-Pierre and Gidon, Albert and Graham, Joe W and
Gupta, Anirudh and Haenel, Valentin and Hay, Etay and Heinis,
Thomas and Hernando, Juan B and Hines, Michael and Kanari, Lida
and Keller, Daniel and Kenyon, John and Khazen, Georges and Kim,
Yihwa and King, James G and Kisvarday, Zoltan and Kumbhar, Pramod
and Lasserre, S{\'e}bastien and Le B{\'e}, Jean-Vincent and
Magalh{\~a}es, Bruno R C and Merch{\'a}n-P{\'e}rez, Angel and
Meystre, Julie and Morrice, Benjamin Roy and Muller, Jeffrey and
Mu{\~n}oz-C{\'e}spedes, Alberto and Muralidhar, Shruti and
Muthurasa, Keerthan and Nachbaur, Daniel and Newton, Taylor H and
Nolte, Max and Ovcharenko, Aleksandr and Palacios, Juan and
Pastor, Luis and Perin, Rodrigo and Ranjan, Rajnish and Riachi,
Imad and Rodr{\'\i}guez, Jos{\'e}-Rodrigo and Riquelme, Juan Luis
and R{\"o}ssert, Christian and Sfyrakis, Konstantinos and Shi,
Ying and Shillcock, Julian C and Silberberg, Gilad and Silva,
Ricardo and Tauheed, Farhan and Telefont, Martin and
Toledo-Rodriguez, Maria and Tr{\"a}nkler, Thomas and Van Geit,
Werner and D{\'\i}az, Jafet Villafranca and Walker, Richard and
Wang, Yun and Zaninetta, Stefano M and DeFelipe, Javier and Hill,
Sean L and Segev, Idan and Sch{\"u}rmann, Felix},
}

@article{nnneighbourhoods,
    author = {Conceição, Pedro and Govc, Dejan and Lazovskis, Jānis and Levi, Ran and Riihimäki, Henri and Smith, Jason P.},
    title = {An application of neighbourhoods in digraphs to the classification of binary dynamics},
    journal = {Network Neuroscience},
    volume = {6},
    number = {2},
    pages = {528-551},
    year = {2022},
}

@InProceedings{bastianautoencoders,
  title = 	 {Topological Autoencoders},
  author =       {Moor, Michael and Horn, Max and Rieck, Bastian and Borgwardt, Karsten},
  booktitle = 	 {Proceedings of the 37th International Conference on Machine Learning},
  pages = 	 {7045--7054},
  year = 	 {2020},
  editor = 	 {III, Hal Daumé and Singh, Aarti},
  volume = 	 {119},
  series = 	 {Proceedings of Machine Learning Research},
  publisher =    {PMLR},
}

@InProceedings{time_varying_merge_trees,
author="Oesterling, Patrick
and Heine, Christian
and Weber, Gunther H.
and Morozov, Dmitriy
and Scheuermann, Gerik",
editor="Carr, Hamish
and Garth, Christoph
and Weinkauf, Tino",
title="Computing and Visualizing Time-Varying Merge Trees for High-Dimensional Data",
booktitle="Topological Methods in Data Analysis and Visualization IV",
year="2017",
publisher="Springer International Publishing",
}

@inproceedings{time_varying_reeb_graphs,
author = {Edelsbrunner, Herbert and Harer, John and Mascarenhas, Ajith and Pascucci, Valerio},
title = {Time-varying reeb graphs for continuous space-time data},
year = {2004},
isbn = {1581138857},
publisher = {Association for Computing Machinery},
booktitle = {Proceedings of the Twentieth Annual Symposium on Computational Geometry},
pages = {366–372},
numpages = {7},
series = {SCG '04}
}

@Article{sliding_windows,
author={Perea, Jose A.
and Harer, John},
title={Sliding Windows and Persistence: An Application of Topological Methods to Signal Analysis},
journal={Foundations of Computational Mathematics},
year={2015},
volume={15},
number={3},
pages={799-838},
}

@article{flagser,
	AUTHOR = {Lütgehetmann, Daniel and Govc, Dejan and Smith, Jason P. and Levi, Ran},
	TITLE = {Computing Persistent Homology of Directed Flag Complexes},
	JOURNAL = {Algorithms},
	VOLUME = {13},
	YEAR = {2020},
	NUMBER = {1},
	ARTICLE-NUMBER = {19},
}

@misc{flagsercount,
  author = {Jason P. Smith},
  title = {\textsc{Flagser-Count}},
  year = {2024},
  publisher = {GitHub},
  journal = {GitHub repository},
  howpublished = {https://github.com/JasonPSmith/flagser-count},
  commit = {c05d9ae3c277d9c8b095106e662f63a7a7af8eb2}
}

@article{reliability,
	title = {Heterogeneous and higher-order cortical connectivity undergirds efficient, robust, and reliable neural codes},
	journal = {iScience},
	volume = {28},
	number = {1},
	pages = {111585},
	year = {2025},
	author = {Daniela {Egas Santander} and Christoph Pokorny and András Ecker and Jānis Lazovskis 	and Matteo Santoro and Jason P. Smith and Kathryn Hess and Ran Levi and Michael W. Reimann},
}

@Article{morphomics,
author={Colombo, Gloria
and Cubero, Ryan John A.
and Kanari, Lida
and Venturino, Alessandro
and Schulz, Rouven
and Scolamiero, Martina
and Agerberg, Jens
and Mathys, Hansruedi
and Tsai, Li-Huei
and Chach{\'o}lski, Wojciech
and Hess, Kathryn
and Siegert, Sandra},
title={A tool for mapping microglial morphology, morphOMICs, reveals brain-region and sex-dependent phenotypes},
journal={Nature Neuroscience},
year={2022},
volume={25},
number={10},
pages={1379-1393},
}

@article{numpy,
 title         = {Array programming with {NumPy}},
 author        = {Charles R. Harris and K. Jarrod Millman and St{\'{e}}fan J.
                 van der Walt and Ralf Gommers and Pauli Virtanen and David
                 Cournapeau and Eric Wieser and Julian Taylor and Sebastian
                 Berg and Nathaniel J. Smith and Robert Kern and Matti Picus
                 and Stephan Hoyer and Marten H. van Kerkwijk and Matthew
                 Brett and Allan Haldane and Jaime Fern{\'{a}}ndez del
                 R{\'{i}}o and Mark Wiebe and Pearu Peterson and Pierre
                 G{\'{e}}rard-Marchant and Kevin Sheppard and Tyler Reddy and
                 Warren Weckesser and Hameer Abbasi and Christoph Gohlke and
                 Travis E. Oliphant},
 year          = {2020},
 journal       = {Nature},
 volume        = {585},
 number        = {7825},
 pages         = {357--362},
}
